\documentclass[11 pt]{article}
\title{Cutting and Pasting in the Torelli Group}
\author{Andrew Putman}
\usepackage{amsmath}
\usepackage{amsfonts}
\usepackage{amssymb}
\usepackage{epsfig}
\usepackage{geometry}
\usepackage{pinlabel}

\newtheorem{theorem}{Theorem}[section]
\newtheorem{theoremsummary}[theorem]{Theorem Summary}

\newtheorem{lemma}[theorem]{Lemma}
\newtheorem{corollary}[theorem]{Corollary}

\newenvironment{proofof}[1]{{\bf {\smallskip}{\noindent}Proof of #1: }}{$\square$ \smallskip}
\newenvironment{proof}{{\bf {\smallskip}{\noindent}Proof: }}{$\square$ \smallskip}
\newenvironment{definition}{{\bf {\smallskip}{\noindent}Definition: }}{\smallskip}

\newenvironment{remark}{{\bf {\smallskip}{\noindent}Remark: }}{\smallskip}
\newenvironment{outline}{{\bf {\smallskip}{\noindent}Outline and Conventions: }}{\smallskip}
\newenvironment{acknowledgements}{{\bf {\smallskip}{\noindent}Acknowledgments: }}{\smallskip}
\newenvironment{history}{{\bf {\smallskip}{\noindent}History and Comments: }}{\smallskip}

\newenvironment{claim}[1]{{\bf {\smallskip}{\noindent}Claim #1: }}{}
\newenvironment{claimproof}{{\bf {\noindent}Proof of Claim: }}{$\square$}

\newcommand\Ker{\mbox{ker}}

\newcommand\Mod{\mbox{Mod}}
\newcommand\Torelli{\mbox{${\mathcal I}$}}

\newcommand\Sp{\mbox{Sp}}

\newcommand\HBolic{\mbox{$\mathbb{H}$}}

\newcommand\Z{\mbox{$\mathbb{Z}$}}
\newcommand\Q{\mbox{$\mathbb{Q}$}}

\newcommand\HH{\mbox{H}}

\newcommand\Max{\mbox{max}}

\newcommand\CaptionSpace{\hspace{0.2in}}
\newcommand\co{\colon\thinspace}

\newcommand\Figure[3]{
\begin{figure}[t]
\centerline{\psfig{file=#2,scale=58}}
\caption{#3}
\label{#1}
\end{figure}}

\newcommand\PtPsh[1]{\mbox{$\sigma_{#1}$}}
\newcommand\LPtPsh[1]{\mbox{$\tilde{\sigma}_{#1}$}}
\newcommand\Curves{\mbox{$\mathcal{C}$}}
\newcommand\CNoSep{\mbox{$\mathcal{C}_{\text{nosep}}$}}
\newcommand\Span[1]{\mbox{$\langle \mbox{#1} \rangle$}}
\newcommand\ComplexLines{\mbox{$\mathcal{L}$}}
\newcommand\AugComplexLines{\mbox{$\overline{\mathcal{L}}$}}
\newcommand\Link{\mbox{link}}
\newcommand\Rk{\mbox{rk}}

\newcommand\Sur{\mbox{Sur}}
\newcommand\TSur{\mbox{TSur}}
\newcommand\Bdry{\mbox{$\partial$}}
\newcommand\Interior{\mbox{Int}}

\begin{document}

\maketitle

\begin{abstract}
We introduce machinery to allow ``cut-and-paste''-style inductive
arguments in the Torelli subgroup of the mapping class group.  In
the past these arguments have been problematic because restricting
the Torelli group to subsurfaces gives different groups depending
on how the subsurfaces are embedded.  We define a category $\TSur$ whose
objects are surfaces together with a decoration restricting how they can be embedded
into larger surfaces and whose morphisms are embeddings which respect the
decoration.  There is a natural ``Torelli functor'' on this category which
extends the usual definition of the Torelli group on a closed surface.  Additionally,
we prove an analogue of the Birman exact sequence for the Torelli groups of
surfaces with boundary and use the action of the Torelli group on the complex
of curves to find generators for the Torelli group.  For genus $g \geq 1$ only
twists about (certain) separating curves and bounding pairs are needed, while
for genus $g=0$ a new type of generator (a ``commutator of a simply intersecting
pair'') is needed.  As a special case, our methods provide a new,
more conceptual proof of the classical result
of Birman-Powell which says that the Torelli group on a closed surface
is generated by twists about separating curves and bounding pairs.
\end{abstract}

\section{Introduction}
Let $\Sigma_{g,n}$ be a genus $g$ surface with $n$ boundary components (we will often omit the $n$ if it equals $0$ and
omit both subscripts if they are unimportant).  The
{\em mapping class group} $\Mod(\Sigma_{g,n})$ is the
group of orientation-preserving homeomorphisms of $\Sigma_{g,n}$ which fix the boundary pointwise
modulo isotopies fixing the boundary pointwise.  The action of
$\Mod(\Sigma_{g,n})$ on $\HH_1(\Sigma_{g,n};\Z)$ preserves the algebraic intersection form.
If $n \leq 1$, then this is a nondegenerate alternating form, so in this case the action yields a
representation $\Mod(\Sigma_{g,n}) \rightarrow \Sp(2g,\Z)$, which is well-known to be surjective.
Its kernel is the {\em Torelli group} $\Torelli(\Sigma_{g,n})$.  Summarizing, for $n \leq 1$ we have an exact sequence
$$1 \longrightarrow \Torelli(\Sigma_{g,n}) \longrightarrow \Mod(\Sigma_{g,n}) \longrightarrow \Sp(2g,\Z) \longrightarrow 1.$$
The group $\Torelli(\Sigma_{g,n})$ plays an important role in both algebraic geometry and low
dimensional topology.  For a survey of the Torelli group (especially the remarkable results of Johnson
which appear in \cite{J2,J3,J4}), see \cite{J1}.

If $i \co \Sigma \hookrightarrow \Sigma'$ is an embedding, then there
is an induced map $i_{\ast} \co \Mod(\Sigma) \rightarrow \Mod(\Sigma')$.  Namely,
if $f \in \Mod(\Sigma)$, then $i_\ast (f)$ equals $f$ on $\Sigma \subset \Sigma'$ and the identity
elsewhere.  In fancier language, if $\Sur$ is the category whose objects are surfaces and whose
morphisms are embeddings, then $\Mod$ is a functor from $\Sur$ to the category of groups and homomorphisms. These morphisms
are fundamental tools in the study of $\Mod(\Sigma)$ (they allow proofs by ``cutting and inducting'').
In this paper, we develop such technology for the Torelli group.

This technology has been problematic in the past because the naive definition of the Torelli
group on a surface with boundary simply does not work.  Indeed, no single definition of $\Torelli(\Sigma_{g,n})$
for $n>1$ satisfies the following two properties, which are the minimum needed for inductive proofs:
\begin{itemize}
\item $\Torelli$ should be functorial in the sense that if $i \co \Sigma \hookrightarrow \Sigma'$ is an embedding,
then $i_{\ast}(\Torelli(\Sigma)) \subset \Torelli(\Sigma')$.
\item $\Torelli$ should be closed under restriction in the sense that if $i \co \Sigma \hookrightarrow \Sigma'$ is
an embedding, then $\Torelli(\Sigma) = i_{\ast}^{-1} (\Torelli(\Sigma'))$.
\end{itemize}
To see that these properties are mutually contradictory and to explain our solution, we need some more concepts.
For a simple closed curve $\gamma$, let $T_{\gamma}$ be the right Dehn twist about $\gamma$.  A curve $\gamma$ is a {\em separating curve} if it separates
the surface into two pieces (for instance, the curve $\gamma_1$ in Figure \ref{figure:differenttorelli}.c).  
A pair of disjoint non-isotopic simple closed curves $\{\gamma,\gamma'\}$ form a {\em bounding
pair} if neither $\gamma$ nor $\gamma'$ separate the surface but $\gamma \cup \gamma'$ does (for instance,
the pair $\{\gamma_2,\gamma_3\}$ in Figure \ref{figure:differenttorelli}.c).  It is not hard to see that
$T_{\gamma} \in \Torelli(\Sigma_g)$ if and only if $\gamma$ is a separating curve, and
similarly if $\{\gamma,\gamma'\}$ is a bounding pair in $\Sigma_{g}$, then $T_{\gamma} T_{\gamma'}^{-1} \in \Torelli(\Sigma_{g})$
(such a mapping class will be called a {\em twist about a bounding pair}).

\Figure{figure:differenttorelli}{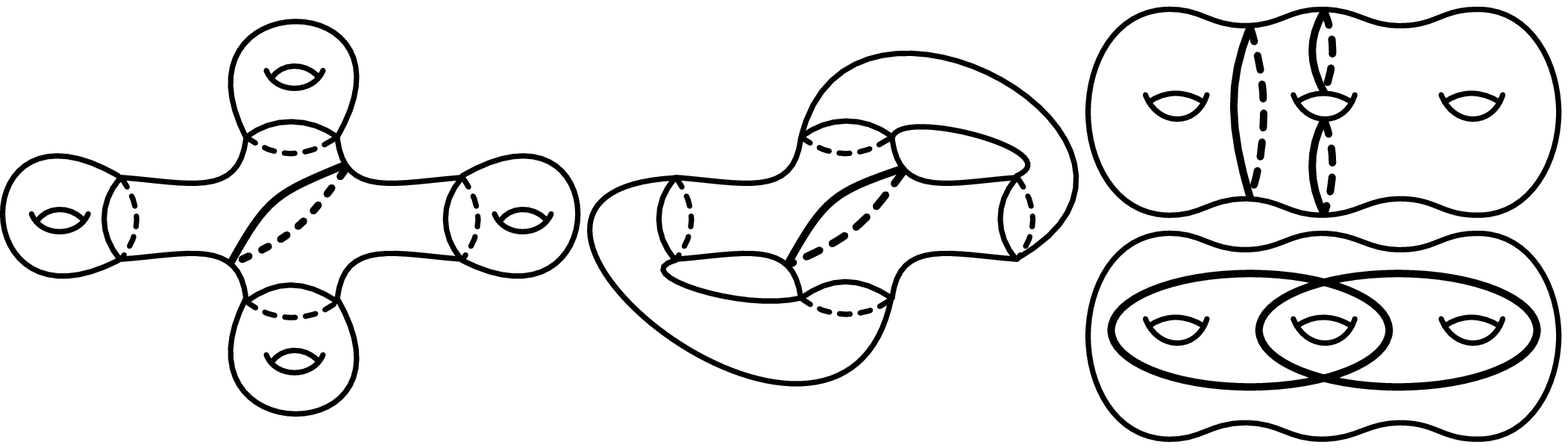}{a,b. Two different embeddings of $\Sigma_{0,4}$ into larger surfaces \CaptionSpace 
c. A separating curve, a bounding pair, and a simply intersecting pair}

Now assume that we have a definition of the Torelli group which is both functorial and closed under restriction, and 
consider Figure \ref{figure:differenttorelli}, which shows two embeddings of $\Sigma_{0,4}$ into closed surfaces.  
The twist $T_{\gamma}$ is a twist about a separating curve in $\Torelli(\Sigma_4)$, so since $\Torelli$ is closed
under restriction, we must have $T_{\gamma} \in \Torelli(\Sigma_{0,4})$.  However, by functoriality we would then
have $T_{\gamma} \in \Torelli(\Sigma_{2})$.  Since $\gamma$ is not a separating curve in $\Sigma_2$, this is a
contradiction.

One difference between the embeddings $\Sigma_{0,4} \hookrightarrow \Sigma_2$ and $\Sigma_{0,4} \hookrightarrow \Sigma_4$ depicted in 
Figure \ref{figure:differenttorelli} is that the partitions
$$\{\partial S \text{ $|$ $S$ is a component of $\Sigma_{i} \setminus \Interior(\Sigma_{0,4})$}\}$$
of the boundary components of $\Sigma_{0,4}$ are different.  It turns out that this additional piece of data
is exactly what we need to distinguish between the behavior of the Torelli groups under different embeddings.  We will
call a pair $(\Sigma,P)$ consisting of a surface $\Sigma$ and a partition $P$ of the boundary components of $\Sigma$
a {\em partitioned surface} (we think of each partition element as a ``chunk'' of the boundary to which we can
attach pieces).  We will construct a category $\TSur$ whose elements are partitioned surfaces and
whose morphisms are embeddings which ``respect the partitions'' (see Section \ref{section:torellimorphisms} for the
precise definition of the morphisms).  We will also define on a partitioned surface $(\Sigma,P)$ a ``homology group''
$\HH_1^P(\Sigma)$ (see Section \ref{section:definition}).  Our main theorem will then
be the following:

\begin{theoremsummary}
\label{theorem:theoremsummary}
There is a functor $\Torelli$ from the category $\TSur$ to the category of groups
and homomorphisms (see Corollary \ref{corollary:torellifunctor}) which satisfies the following properties:
\begin{itemize}
\item If $i \co \Sigma \hookrightarrow \Sigma_g$ is an embedding and
$$P = \{\partial S  \text{ $|$ $S$ is a component of $\Sigma_{g} \setminus \Interior(\Sigma)$}\},$$
then $\Torelli(\Sigma,P) = i_{\ast}^{-1}(\Torelli(\Sigma_g))$ (see Theorem \ref{theorem:torelliwelldefined}).
\item For a partitioned surface $(\Sigma,P)$, the group $\Torelli(\Sigma,P)$ equals
the subgroup of $\Mod(\Sigma)$ acting trivially on $\HH_1^P(\Sigma)$ (see Theorem \ref{theorem:torelliwelldefined}).
\item $\Torelli$ is closed under restriction in the following sense : if $(\Sigma',P')$ is a partitioned
surface and $i \co \Sigma \hookrightarrow \Sigma'$ is an embedding, then there is some partition $P$ of the
boundary components of $\Sigma$ so that $\Torelli(\Sigma,P) = i_{\ast}^{-1}(\Torelli(\Sigma',P'))$ (see
Theorem \ref{theorem:torellirestriction}).
\end{itemize}
\end{theoremsummary}

Now, the mapping class groups on surfaces with differing numbers of boundary components
are related by the fundamental Birman exact sequences (see \cite{Bi1, Bi3} and Section \ref{section:birmanexactsequence}).  
One version of this, due to Johnson \cite{J2}, is that for $n \geq 1$ and $(g,n) \neq (1,1)$, we have an exact sequence
$$1 \longrightarrow \pi_1(U\Sigma_{g,n-1}) \longrightarrow \Mod(\Sigma_{g,n}) \longrightarrow \Mod(\Sigma_{g,n-1}) \longrightarrow 1,$$
where $U\Sigma_{g,n-1}$ is the unit tangent bundle of $\Sigma_{g,n-1}$, the map $\Mod(\Sigma_{g,n}) \longrightarrow \Mod(\Sigma_{g,n-1})$
is induced by ``gluing a disc to a boundary component'', and the map $\pi_1(U\Sigma_{g,n-1}) \rightarrow \Mod(\Sigma_{g,n})$ is
induced by ``pushing the boundary component around curves''.  We construct a similar exact sequence for the Torelli
group:
\begin{theorem}
\label{theorem:torelliexactsequenceprelim}
Let $(\Sigma_{g,n},P)$ be a partitioned surface with $n \geq 1$ and $(g,n) \neq (1,1)$, and let $b$ be a boundary component of $\Sigma_{g,n}$ contained
in a set $p \in P$.  Consider the map $i \co \Sigma_{g,n} \hookrightarrow \Sigma_{g,n-1}$ induced by gluing a disc to $b$, and let $P'$ be the partition
of the boundary components of $\Sigma_{g.n-1}$ induced by $i$.  There is then an exact sequence
$$1 \longrightarrow K \longrightarrow \Torelli(\Sigma_{g,n},P) \longrightarrow \Torelli(\Sigma_{g,n-1},P') \longrightarrow 1$$
with $K$ equal to the following:
\begin{itemize}
\item If $p = \{b\}$, then $K = \pi_1(U\Sigma_{g,n-1})$.
\item If $p \neq \{b\}$, then $K$ is isomorphic to the kernel of the natural map $\pi_1(\Sigma_{g,n-1}) \rightarrow \HH_1^{P'}(\Sigma_{g,n-1})$ (recall
that $\HH_1^{P'}(\Sigma_{g,n-1})$ is the ``homology group'' discussed in Theorem \ref{theorem:theoremsummary}).
\end{itemize}
\end{theorem}
See Section \ref{section:torellisequence} for a discussion of the manner in which the kernel $K$ embeds into
$\pi_1(U\Sigma_{g,n-1})$.

Finally, we combine the machinery we have developed with action of the Torelli group on the complex
of curves $\Curves(\Sigma)$ (see Section \ref{section:torelligenerators}) to determine generators for $\Torelli(\Sigma,P)$.
We define ``$P$-separating curves'' and ``$P$-bounding pairs'' to be separating curves and bounding pairs
which ``respect the partition $P$'' (see the end of Section \ref{section:definition} for a precise
definition).  Our first theorem is then the following: 
\begin{theorem}
\label{theorem:torelligenerators}
For any partitioned surface $(\Sigma_{g,n},P)$ with $g \geq 1$, the group $\Torelli(\Sigma_{g,n},P)$ is generated
by twists about $P$-separating curves and twists about $P$-bounding pairs.
\end{theorem}
As a special case of this theorem, we obtain a new proof of the following classical
theorem of Birman-Powell \cite{Bi2,P}:
\begin{corollary}
For a surface $\Sigma_{g,n}$ with $n \leq 1$, the group $\Torelli(\Sigma_{g,n})$ is generated
by twists about separating curves and bounding pairs.
\end{corollary}
Our proof of Theorem \ref{theorem:torelligenerators} makes use of a basic result of Armstrong \cite{A} (see
Section \ref{section:armstrong}) which says that if a group $G$ acts nicely on a simply-connected simplicial
complex $X$, then $G$ is generated by elements which stabilize vertices of $X$ if and only if $X/G$ is simply
connected.  We prove that $\Curves(\Sigma_{g}) / \Torelli(\Sigma_{g})$ is simply connected for $g \geq 2$, so
Armstrong's theorem says that $\Torelli(\Sigma_{g})$ is generated by the stabilizer subgroups $(\Torelli(\Sigma_{g}))_{\gamma}$ of
simple closed curves $\gamma$.  These stabilizer subgroups are supported on ``simpler'' subsurfaces, so we can use induction
and Theorem \ref{theorem:torelliexactsequenceprelim} to analyze them.

Now, the condition $g \geq 1$ in Theorem \ref{theorem:torelligenerators} is necessary; indeed, for some partitions
$P$ of the boundary components of $\Sigma_{0,n}$, there are {\em no} $P$-separating curves or $P$-bounding pairs.  To
find generators for the Torelli groups of genus $0$ surfaces, we make one final definition.  If $\{\gamma,\gamma'\}$
is a pair of simple closed curves in $\Sigma_{g,n}$ whose geometric intersection number is $2$ and whose algebraic intersection
number is $0$ (for instance, the pair $\{\gamma_4,\gamma_5\}$ in Figure \ref{figure:differenttorelli}.c), then
it is easy to see that for any partition $P$ of the boundary components of $\Sigma_{g,n}$ we have
$[T_{\gamma},T_{\gamma'}] \in \Torelli(\Sigma_{g,n},P)$.  We will call these {\em commutators of simply intersecting
pairs}.  Our final theorem is the following:

\begin{theorem}
\label{theorem:torelligenerators2}
For any genus $0$ partitioned surface $(\Sigma_{0,n},P)$, the group $\Torelli(\Sigma_{0,n},P)$ is generated
by twists about $P$-separating curves, twists about $P$-bounding pairs, and commutators of simply
intersecting pairs.
\end{theorem}
 
\begin{history}
At least two special cases of our construction have appeared in the literature.  The
simplest appears in the work of Hain \cite{H}, who (in our notation) considered
$\Torelli(\Sigma_{g,n},P)$ with $P=\{\{1\},\ldots,\{n\}\}$, which he defined as the subgroup
of $\Torelli(\Sigma_{g})$ fixing $n$ discs.  This is the ``largest'' possible definition
of the Torelli group, and is rather easily related to the closed surface case.  However,
it does not have good functoriality properties, and it seems difficult to use it
in inductive proofs (though, as Hain notes, it does have interesting interpretations
in terms of algebraic geometry).

The other special case appears in the
work of Johnson \cite{J3} and van den Berg \cite{vdB}.  In our notation, Johnson considered $\Torelli(\Sigma_{g,n},P)$
with $P=\{\{1,\ldots,n\}\}$.  This is the ``smallest'' possible definition of
the Torelli group, and it is functorial under embeddings.  However, it is
not closed under restriction except in the simplest possible cases, and
great care has to be exercised when using it in inductive proofs.  
This work was continued in the unpublished thesis of van den Berg \cite{vdB}, who indicated
how to extend Johnson's calculation of $\HH_1(\Torelli(\Sigma_{g});\Z)$ to surfaces
with boundary using Johnson's definition of the Torelli group on surfaces with boundary.  
She also gave a very brief sketch of the identification
of the kernel of the exact sequence in Theorem \ref{theorem:torelliexactsequenceprelim} 
for the cases she was considering, though she did not prove that the associated map was surjective.

The history of Birman-Powell's result is rather interesting.  Though it is a fundamental result used in
nearly every subsequent paper on the Torelli group, their proof is the only one we are aware of
in the literature.  The story begins with a paper of Klingen \cite{K} in which he gave
an algorithm for computing a presentation of $\Sp_{2g}(\Z)$.  In Birman's
paper \cite{Bi2}, she followed this algorithm.  Assembling results from a paper of Magnus \cite{Ma}
and the unpublished thesis of Gold \cite{G}, she showed that Klingen's presentation has
5 families of generators and 67 families of relations.  After an absolutely
heroic calculation (whose details, needless to say, are only briefly sketched in Birman's
short paper), she reduced this to a presentation with 3 families of generators and 10 families of relations.

These 3 families of generators are the images of a standard set of generators for $\Mod(\Sigma_{g})$.  She
calculated that all but a few of these families of relations lift to relations in $\Mod(\Sigma_{g})$.
The relations in $\Sp_{2g}(\Z)$ which do not lift to
relations in $\Mod(\Sigma_{g})$ lift to normal generators for $\Torelli(\Sigma_{g})$.  In Powell's paper,
he showed how to express these normal generators for $\Torelli(\Sigma_{g})$ as products of twists about
separating curves and bounding pairs, thus establishing the result.

Our proof seems to be the first to appear in the literature which does not depend on
unpublished results and for which full details are given.  Of course, $K$-theoretic methods
have by now yielded simpler presentations of $\Sp_{2n}(\Z)$ than Klingen's
presentation (see, e.g., \cite[Theorem 9.2.13]{HOM}), but these constructions
are quite involved, and it is nontrivial to perform Birman and Powell's analysis
on them to get generators for the Torelli group.  Our method has two advantages
over such an approach.  First, our methods are ``intrinsic'' to the theory of mapping class groups.
Second, for the most part we avoid complicated group-theoretic calculations.
\end{history}

\begin{acknowledgements}
I wish to thank my advisor Benson Farb for introducing me to the Torelli group and commenting
extensively on previous incarnations of this paper,
Ben Wieland for several useful discussions, and Justin Malestein for many useful discussions
and corrections (in particular, Justin observed Lemma \ref{lemma:augcomplexfarey}, which
significantly simplified my original argument).
Additionally, I wish to thank Dan Margalit and Matt Day for some corrections and Yair Minsky for allowing me to use
his drawing of the Farey graph.  Finally, I wish to thank the mathematics
department of the Georgia Institute of Technology for their hospitality during
the time in which part of this paper was conceived.
\end{acknowledgements}

\begin{outline}
We begin with two sections outlining preliminaries.  Next, we define the Torelli group on a surface with boundary
in Section \ref{section:definition} and discuss morphisms between different Torelli groups
in Section \ref{section:torellimorphisms}.  After a discussion of how to restrict Torelli groups
to subsurfaces in Section \ref{subsection:restriction}, we prove our analogue of the Birman
exact sequence in Section \ref{section:torellisequence}.  In Section \ref{section:torelligenerators}
we find generators for the Torelli groups, proving Theorems \ref{theorem:torelligenerators} and
\ref{theorem:torelligenerators2}.  The proofs of these theorems depend on the simple-connectedness of $\Curves(\Sigma_{g}) / \Torelli(\Sigma_{g})$
for $g \geq 2$, which is proven in Section \ref{section:curvestorelli} (see the beginning of that section for an outline
of this lengthy proof).  We conclude with an Appendix in which we prove several useful lemmas
about the topology of surfaces.

Throughout this paper, all homology groups will have $\Z$-coefficients.  
A {\em summand} $A$ of a module $M$ is a submodule so that there
exists another submodule $B$ with $M=A \oplus B$.  We will often abuse
notation and discuss the span $\Span{$X_1, \ldots ,X_k$}$ of a set of submodules 
of a single module $X$.  The algebraic intersection of two homology classes $h_1,h_2$ will
be denoted by $i_a(h_1,h_2)$.  The geometric intersection number of two simple closed curves
$\gamma_1,\gamma_2$ will be denoted by $i_g(\gamma_1,\gamma_2)$.  Finally, the flag complex with vertices $X$ and adjacency relation $R$ is
the simplicial complex whose simplices are sets $\{x_1,\ldots,x_k\} \subset X$ so that $x_i R x_j$ for all $i$ and $j$.
\end{outline}

\section{Preliminaries}
\label{section:preliminaries}

\subsection{The Birman exact sequences}
\label{section:birmanexactsequence}
In this section, we will give a detailed review of the exact sequences of Birman and
Johnson \cite{Bi1,Bi3,J2} which describe the effect on the mapping class group of gluing a
disc to a boundary component.  We will need the following definition:

\begin{definition}
Consider a surface $\Sigma_{g,n}$.  Let $x \in \Sigma_{g,n}$ be a point.  We define $\Mod(\Sigma_{g,n},x)$, the {\em mapping class
group relative to $x$}, to be the group of orientation-preserving homeomorphisms which fix $x$ and
the boundary modulo isotopies fixing $x$ and the boundary.
\end{definition}

\noindent
Let $b$ be a boundary component of $\Sigma_{g,n}$.  There is a natural embedding $\Sigma_{g,n} \hookrightarrow \Sigma_{g,n-1}$
induced by gluing a disc to $b$.  Let $x$ be a point in the interior of the new disc.  Clearly we can
factor the induced map $\Mod(\Sigma_{g,n}) \rightarrow \Mod(\Sigma_{g,n-1})$ into a composition
$$\Mod(\Sigma_{g,n}) \longrightarrow \Mod(\Sigma_{g,n-1},x) \longrightarrow \Mod(\Sigma_{g,n-1}).$$
Now let $U\Sigma_{g,n-1}$ be the unit tangent bundle of $\Sigma_{g,n-1}$ and $\tilde{x}$
be any lift of $x$ to $U\Sigma_{g,n-1}$.  The combined work of Birman \cite{Bi1} and Johnson \cite{J2}
shows that (except for the degenerate case $(g,n)=(1,1)$) all of our groups fit into the following commutative diagram with exact rows and columns:

\begin{center}
\begin{tabular}{c@{\hspace{0.05 in}}c@{\hspace{0.05 in}}c@{\hspace{0.05 in}}c@{\hspace{0.05 in}}c@{\hspace{0.05 in}}c@{\hspace{0.05 in}}c@{\hspace{0.05 in}}c@{\hspace{0.05 in}}c}
    &               & $1$                        &               & $1$                      &               &                        &               &     \\
    &               & $\downarrow$               &               & $\downarrow$             &               &                        &               &     \\
    &               & $\Z$                       & $=$           & $\Z$                     &               &                        &               &     \\
    &               & $\downarrow$               &               & $\downarrow$             &               &                        &               &     \\
$1$ & $\longrightarrow$ & $\pi_1(U\Sigma_{g,n-1},\tilde{x})$ & $\longrightarrow$ & $\Mod(\Sigma_{g,n})$     & $\longrightarrow$ & $\Mod(\Sigma_{g,n-1})$ & $\longrightarrow$ & $1$ \\
    &               & $\downarrow$               &               & $\downarrow$             &               & $\parallel$            &               &     \\
$1$ & $\longrightarrow$ & $\pi_1(\Sigma_{g,n-1},x)$  & $\longrightarrow$ & $\Mod(\Sigma_{g,n-1},x)$ & $\longrightarrow$ & $\Mod(\Sigma_{g,n-1})$ & $\longrightarrow$ & $1$ \\
    &               & $\downarrow$               &               & $\downarrow$             &               &                        &               &     \\
    &               & $1$                        &               & $1$                      &               &                        &               &     
\end{tabular}
\end{center}

\Figure{figure:birmanexactsequence}{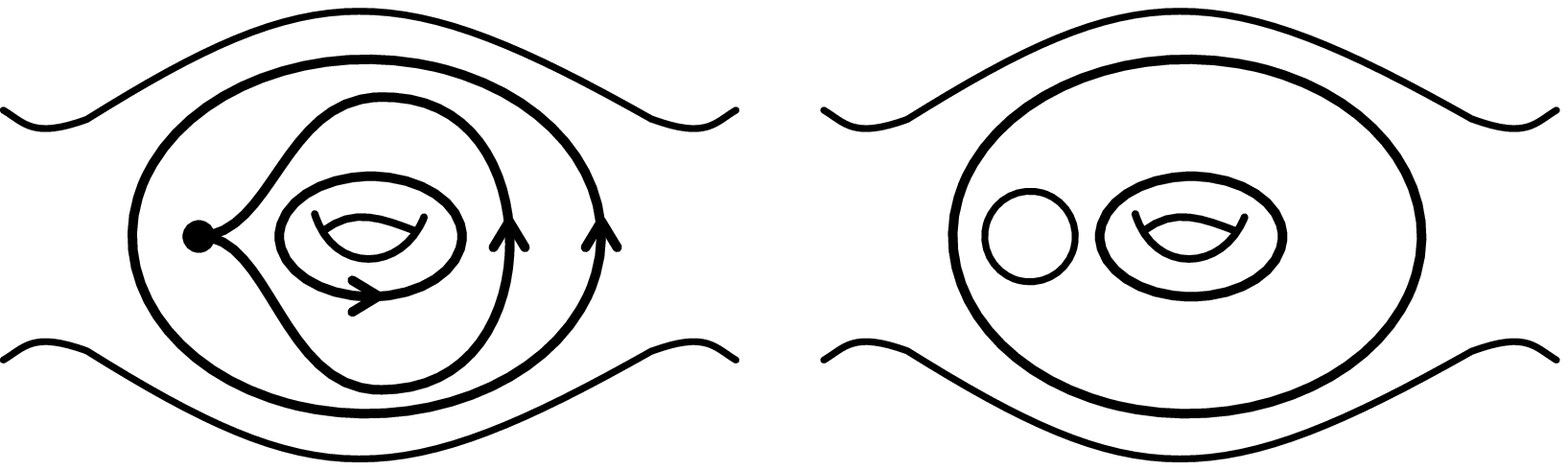}{a. Image of a simple closed curve in $\Mod(\Sigma_{g,n-1},x)$ \CaptionSpace 
b. Lift of a simple closed curve to $\Mod(\Sigma_{g,n})$}

\noindent
The $\Z$ in the first column is the loop in the fiber, while the $\Z$ in the second column corresponds to the Dehn twist about the
filled-in boundary component.  

For $\gamma \in \pi_1(\Sigma_{g,n-1},x)$, let $\PtPsh{\gamma}$ be the element of $\Mod(\Sigma_{g,n-1},x)$
associated to $\gamma$ (hence $\PtPsh{\gamma}$ is the mapping class which ``pushes $x$ around $\gamma$'').  If $\gamma$ is a simple closed 
curve, then there is a nice formula for $\PtPsh{\gamma}$ (see Figure
\ref{figure:birmanexactsequence}.a).  Namely, let $\gamma_1$ and $\gamma_2$ be the boundary components of
a regular neighborhood of $\gamma$.  The orientation of $\gamma$ induces orientations on $\gamma_1$ and $\gamma_2$;
assume that $\gamma$ lies to the left of $\gamma_1$ and to the right of $\gamma_2$.  Then $\PtPsh{\gamma}=T_{\gamma_1} T_{\gamma_2}^{-1}$.

Continue to assume that $\gamma$ is a simple closed curve.  We will construct a natural lift $\LPtPsh{\gamma}$ of
$\PtPsh{\gamma}$ to $\pi_1(U\Sigma_{g,n-1},\tilde{x}) \subset \Mod(\Sigma_{g,n})$ (we reiterate that our construction depends on the simplicity of
$\gamma$).  Recall that we have been
considering $\Sigma_{g,n-1}$ to be $\Sigma_{g,n}$ with a disc glued to $b$.  In the other direction, we can consider
$\Sigma_{g,n}$ to be $\Sigma_{g,n-1}$ with the point $x$ blown up to a boundary component.  Two such identifications
may differ by a power of $T_b$; however, since $T_b$ fixes both $\gamma_1$ and $\gamma_2$, there are well-defined
lifts $\tilde{\gamma}_1$ and $\tilde{\gamma}_2$ of the $\gamma_i$ to $\Sigma_{g,n}$ (see Figure
\ref{figure:birmanexactsequence}.b).  It is not hard to see that
$$\LPtPsh{\gamma} := T_{\tilde{\gamma}_1} T_{\tilde{\gamma}_2}^{-1}$$
is a lift of $\PtPsh{\gamma}$.

\begin{remark}
While the map
$$\PtPsh{\cdot} \co \pi_1(\Sigma_{g,n-1},x) \longrightarrow \Mod(\Sigma_{g,n-1},x)$$
is a homomorphism, the map $\LPtPsh{\cdot}$ (which is only defined on simple closed curves)
does not extend to a homomorphism.
\end{remark}

\subsection{Groups acting on simplicial complexes}
\label{section:armstrong}
In this section, we will prove a theorem (really, a corollary of a theorem of
Armstrong \cite{A}) which we will use to find generating sets for the Torelli group.
We will need the following definition:

\begin{definition}
A group $G$ acts on a simplicial complex $X$ {\em without rotations} if for all
simplices $s$ of $X$ the stabilizer $G_s$ stabilizes $s$ pointwise.
\end{definition}

\noindent
The following is our theorem:

\begin{theorem}
\label{theorem:armstrong}
Let $G$ act without rotations on a simply-connected simplicial complex $X$.  Then
$G$ is generated by the set
$$\bigcup_{v \in X^{(0)}} G_v$$
if and only if $X/G$ is simply-connected.
\end{theorem}
\begin{proof}
Let $H$ be the subgroup of $G$ generated by the indicated set (note that $H$ is normal).  Observe that since
$G$ acts without rotations, we can subdivide $X$ without changing
$H$.  After possibly subdividing twice, the work of Armstrong \cite{A} implies that there is an exact sequence
$$1 \longrightarrow H \longrightarrow G \longrightarrow \pi_1(X/G) \longrightarrow 1.$$
This clearly implies the theorem.  We briefly indicate the construction of the exact sequence.  Let
$\pi \co X \rightarrow X/G$ be the projection, and fix
a vertex $v$ of $X$.  We define a homomorphism $j \co G \rightarrow \pi_1(X/G,\pi(v))$ in the following
way.  For $g \in G$, let $\ell$ be a path in $X$ from $v$ to $gv$.  Then $j(g):=\pi \circ \ell$ is
a loop in $X/G$ based at $\pi(v)$.  Since $X$ is simply-connected, $j$ is well-defined, and
it is clear that $j$ is a surjective homomorphism.  Now, if
$gw=w$ for a vertex $w \in X$, we claim that $j(g)=1$.  Indeed, let $\ell'$ be a path from $v$
to $w$.  Then $\ell' \cdot (g (\ell'))^{-1}$ is a path from $v$ to $gv$ which clearly projects to a
null-homotopic loop in $X/G$, proving the claim.  The bulk of Armstrong's work, therefore, consists
of showing that if $g \in \Ker(j)$, then $g \in H$.  We refer the reader to Armstrong's
paper \cite[pp. 643-645]{A} for the details.
\end{proof}

\section{Definition of the Torelli group on a surface with boundary}
\subsection{Definition}
\label{section:definition}
In this section, we will define the Torelli group on a surface with boundary.  As discussed
in the introduction, our main goal is to understand the subgroups of $\Mod(\Sigma)$ which
arise as $i_{\ast}^{-1}(\Torelli(\Sigma_{g}))$ for embeddings $i \co \Sigma \hookrightarrow \Sigma_{g}$.
We observed there that different embeddings yield different ``Torelli groups'' for $\Sigma$.
 
\Figure{figure:capping}{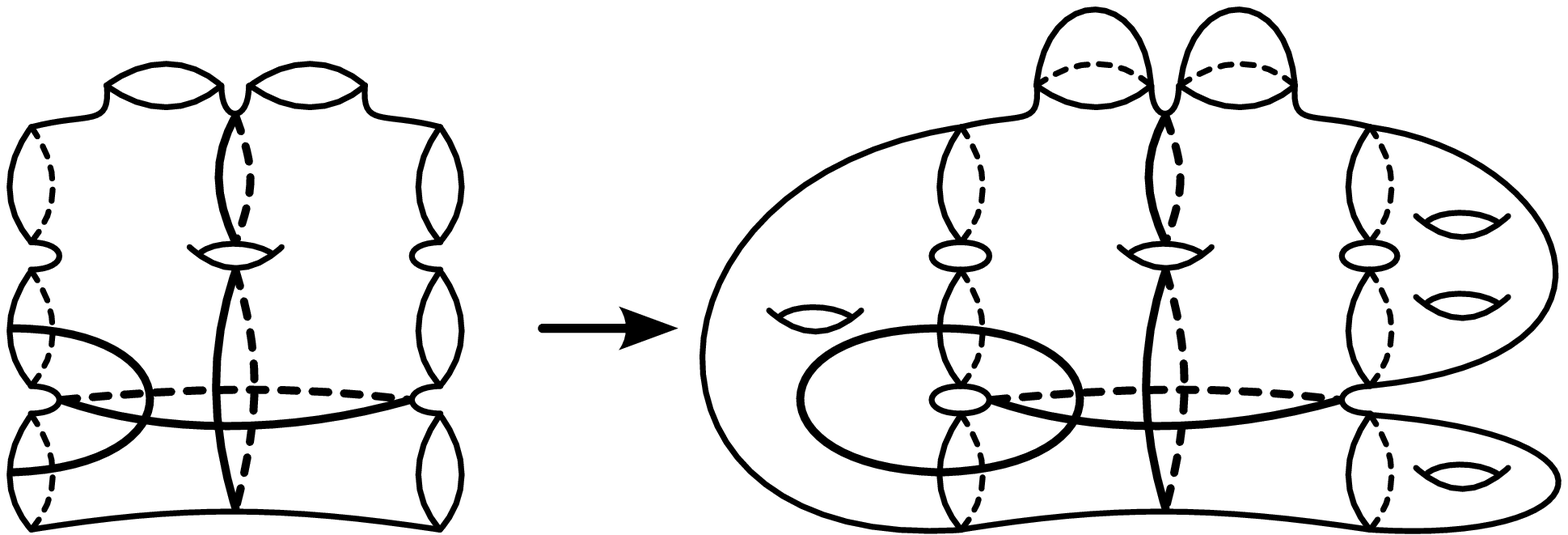}{A capping of $(\Sigma_{1,8},\{\{1,2,3\},\{4,5\},\{6\},\{7\},\{8\}\})$}

Recall that in the introduction we indicated that to distinguish these different Torelli groups
we would attach to a surface $\Sigma$ a partition $P$ of its boundary components; the pair
$(\Sigma,P)$ will be called a {\em partitioned surface}.  Associated to an embedding
$i \co \Sigma \hookrightarrow \Sigma_g$ is a partition
$$P = \{\partial S  \text{ $|$ $S$ is a component of $\Sigma_{g} \setminus \Interior(\Sigma)$}\}$$
of the boundary components of $\Sigma$.  The following is a useful shorthand for the embeddings which give rise to a partition $P$:

\begin{definition}
A {\em capping} of a partitioned surface $(\Sigma,P)$ (see Figure \ref{figure:capping}) is an
embedding
$$\Sigma \hookrightarrow \Sigma_{g}$$
so that for each component $S$ of $\Sigma_{g} \setminus \Interior(\Sigma)$, the
set of boundary components of $S$ is exactly equal to an element of $P$.
\end{definition}

\noindent
The follow obvious lemma says that that every embedding is a capping for an
appropriate partition:

\begin{lemma}
\label{lemma:embeddingsarecappings}
Let $i \co \Sigma \hookrightarrow \Sigma_{g}$ be an embedding.  Set
$$P = \{\partial S  \subset \partial \Sigma \text{ $|$ $S$ is a component of $\Sigma_{g} \setminus \Interior(\Sigma)$}\}.$$
Then $i$ is a capping of $(\Sigma,P)$.
\end{lemma}

\noindent
We now define the Torelli group of a partitioned surface.

\begin{definition}
For a partitioned surface $(\Sigma,P)$, let $\Torelli(\Sigma,P):=i_{\ast}^{-1}(\Torelli(\Sigma_g))$
for any capping $i \co \Sigma \hookrightarrow \Sigma_{g}$.
\end{definition}

\noindent
Of course, it is not at all clear that this definition of $\Torelli(\Sigma,P)$ is independent
of the chosen capping.  Also, it seems rather ad-hoc -- one would like to say that $\Torelli(\Sigma,P)$
is the subgroup of $\Mod(\Sigma)$ fixing some sort of homology group.  We will answer these
objections simultaneously by defining a certain intrinsic ``homology group'' $\HH_1^P(\Sigma)$
and then proving Theorem \ref{theorem:torelliwelldefined}, which says that $\Torelli(\Sigma,P)$ is
exactly the subgroup of $\Mod(\Sigma)$ acting trivially on $\HH_1^P(\Sigma)$.

The construction of $\HH_1^P(\Sigma)$ is a two step process.  For the first step,
observe that in Figure \ref{figure:capping} the mapping class $T_{\gamma_1} T_{\gamma_2}^{-1}$ is an element
of $\Torelli(\Sigma_{8})$, but it does not fix $\HH_1(\Sigma_{1,8})$.  The problem is that (picking
appropriate orientations for everything) 
$$[\gamma_1] - [\gamma_2] = [b_1]+[b_2]+[b_3]+[b_7] \neq 0.$$
This motivates the following definition:

\begin{definition}
Consider a partitioned surface $(\Sigma,P)$, and enumerate the partition $P$:
$$P=\{\{b_1^1,\ldots,b_{k_1}^1\},\ldots,\{b_1^m,\ldots,b_{k_m}^m\}\}.$$
Orient the boundary components $b_i^j$ so that $\sum_{i,j} [b_i^j] = 0$ in $\HH_1(\Sigma)$.  Define
\begin{align*}
\Bdry \HH_1^P(\Sigma) &= \Span{$[b_1^1]+\ldots+[b_{k_1}^1],\ldots,[b_1^m]+\ldots+[b_{k_m}^m]$} \subset \HH_1(\Sigma), \\
\overline{\HH}_1^P(\Sigma) &= \HH_1(\Sigma) / \Bdry \HH_1^P(\Sigma).
\end{align*}
\end{definition}

\noindent
The following lemma is immediate from the definitions:

\begin{lemma}
\label{lemma:cappingh1injective}
If $i \co \Sigma \hookrightarrow \Sigma_{g}$ is a capping of a partitioned surface $(\Sigma,P)$, then there is an induced
injection $i_{\ast} \co \overline{\HH}_1^P(\Sigma) \hookrightarrow \HH_1(\Sigma_{g})$.
\end{lemma}

\noindent
For the second step of the construction of $\HH_1^P(\Sigma)$, observe that in Figure
\ref{figure:capping} the mapping class $T_{\delta}$ does not fix the homology class $[h_1]+[h_2]$ (where the $h_i$
are the indicated arcs and $[h_i]$ is the chain corresponding to $h_i$; only the sum $[h_1]+[h_2]$ is a cycle), and hence
is not in $\Torelli(\Sigma_{8})$.  However, it does fix $\overline{\HH}_1^P(\Sigma)$.  The difficulty is
that we really need elements of $\Torelli(\Sigma_{1,8},\{\{1,2,3\},\{4,5\},\{6\},\{7\},\{8\}\})$ to fix
$[h_1]$, which is a homology class relative to points on the boundary components.  This motivates the
following definition:

\begin{definition}
Fix a partitioned surface $(\Sigma,P)$, and let $Q$ be a set containing one point from each boundary component of $\Sigma$.  Define
$\HH_1^P(\Sigma)$ to equal the image of
the following submodule of $\HH_1(\Sigma,Q)$ in $\HH_1(\Sigma,Q) / \Bdry \HH_1^P(\Sigma)$:
\begin{align*}
\langle \{[h] \in \HH_1(\Sigma,Q) \text{ $|$ } &\text{either $h$ is a simple closed curve or $h$ is a properly} \\
                                               &\text{embedded arc from $q_1$ to $q_2$ with $q_1,q_2 \in Q$ lying} \\
                                               &\text{in boundary components $b_1$ and $b_2$ with $\{b_1,b_2\} \subset p$} \\
                                               &\text{for some $p \in P$\}} \rangle.
\end{align*}
We remark that in the future we will omit mention of $Q$ and instead refer to the homology classes of
arcs between boundary components, the set $Q$ being understood.  If  
boundary components $b_1$ and $b_2$ satisfy $b_1,b_2 \in p$ for some $p \in P$, we will call them 
{\em $P$-adjacent boundary components}.
\end{definition}

\begin{remark}
The idea of using relative homology classes to analyze the Torelli groups on
surfaces with boundary is due to Johnson \cite{J3}.
\end{remark}

It is clear that $\Mod(\Sigma)$ acts upon $\HH_1^P(\Sigma)$.  We now prove the following:

\begin{theorem}
\label{theorem:torelliwelldefined}
For a partitioned surface $(\Sigma,P)$, the group $\Torelli(\Sigma,P)$ is exactly the subgroup of
$\Mod(\Sigma)$ which acts trivially on $\HH_1^P(\Sigma)$.  In particular, $\Torelli(\Sigma,P)$
is well-defined.
\end{theorem}
\begin{proof}
Fix a capping $i \co \Sigma \hookrightarrow \Sigma_{g}$ of $(\Sigma,P)$.  Define the following subsets of $\HH_1(\Sigma_g)$:
\begin{align*}
Q_1 = \{[h] \in \HH_1(\Sigma_{g})          \text{ $|$ } &\text{$h$ is a simple closed curve in $\Sigma_{g} \setminus \Sigma$}\}, \\
Q_2 = \{[h] \in \HH_1(\Sigma_{g})          \text{ $|$ } &\text{$h$ is a simple closed curve in $\Sigma$}\}, \\
Q_3 = \{[h_1]+[h_2] \in \HH_1(\Sigma_{g})  \text{ $|$ } &\text{$h_1$ is a properly embedded arc in $\Sigma$ between} \\
                                                         &\text{$P$-adjacent boundary components and $h_2$ is} \\
                                                         &\text{a properly embedded arc in $\Sigma_{g} \setminus \Interior(\Sigma)$} \\
                                                         &\text{with the same endpoints as $h_1$}\}.
\end{align*}
For an example of an element of $Q_3$, see Figure \ref{figure:capping}.  It is clear that
$$\langle Q_1 \cup Q_2 \cup Q_3 \rangle = \HH_1(\Sigma_{g}).$$
For $f \in \Mod(\Sigma)$, the mapping class $i_\ast (f)$ fixes every element of $Q_1$.  Also, by Lemma
\ref{lemma:cappingh1injective}, the mapping class $i_\ast(f)$ fixes $[h] \in Q_2$ if and only if $f$ fixes the corresponding element
of $\HH_1^P(\Sigma)$.  Finally, we claim that $i_{\ast}(f)$ fixes $[h_1]+[h_2] \in Q_3$ if and only
if $f$ fixes $[h_1] \in \HH_1^P(\Sigma)$.  Indeed, the reverse implication is trivial, while for the forward
implication, observe that if $f$ does not fix $[h_1] \in \HH_1^P(\Sigma)$, then $[h_1]-f([h_1]) \in \overline{\HH}_1^P(\Sigma)$ is nonzero,
so by Lemma \ref{lemma:cappingh1injective} we have 
$$([h_1]+[h_2])-i_{\ast}(f)([h_1]+[h_2])=i([h_1]-f([h_1])) \neq 0;$$
i.e. $i_{\ast}(f)$ does not fix $[h_1]+[h_2]$.  We conclude that
$i_\ast(f)$ acts trivially on $\HH_1(\Sigma_{g})$ if and only if $f$ acts trivially on $\HH_1^P(\Sigma)$, as desired.
\end{proof}

\noindent
We now prove the following lemma, which describes the natural ``algebraic
intersection pairing'' on $\HH_1^P(\Sigma)$:

\begin{lemma}
\label{lemma:homologypairing}
Fix a partitioned surface $(\Sigma,P)$, and let $Q \subset \partial \Sigma$ be the set
from the definition of $\HH_1^P(\Sigma)$.  Then the algebraic intersection pairing $i_a(\cdot,\cdot)$
on $\HH_1(\Sigma,Q)$ induces a pairing on $\HH_1^P(\Sigma)$ (which we will continue to call $i_a(\cdot,\cdot)$).
This pairing is preserved by $\Mod(\Sigma)$, and for a simple closed curve $\gamma$, the mapping
class $T_{\gamma}$ acts upon $\HH_1^P(\Sigma)$ by the transvection
$$h \longmapsto h + i_a(h,[\gamma]) [\gamma].$$
\end{lemma}
\begin{proof}
The only nontrivial part of this lemma is that $i_a$ is well defined on $\HH_1^P(\Sigma)$.  Let
$H \subset \HH_1(\Sigma,Q)$ be the pull-back of $\HH_1^P(\Sigma)$.  We must show that for
$b \in \partial \HH_1^P(\Sigma)$, the map $i_a(b,\cdot)$ restricts to the $0$ map on $H$.  Enumerating the
partition $P$ as
$$P=\{\{b_1^1,\ldots,b_{k_1}^1\},\ldots,\{b_1^m,\ldots,b_{k_m}^m\}\},$$
we can assume that $b$ is one of the generators $[b_1^i]+\cdots+[b_{k_i}^i]$ of $\partial \HH_1^P(\Sigma)$.
If $h$ is a simple closed curve or an arc between boundary
components $b_r^j$ and $b_s^j$ with $i \neq j$, then trivially $i_a(b,[h])=0$.  If instead $h$
is an arc between boundary components $b_r^i$ and $b_s^i$, then
we have $i_a(b,[h]) = i_a([b_r^j],[h]) + i_a([b_s^j],[h])=0$.  Since $i_a(b,\cdot)$ vanishes
on generators for $H$, it vanishes on $H$, as desired.
\end{proof}

\noindent
We conclude this section by discussing typical elements of $\Torelli(\Sigma,P)$.

\begin{definition}
Fix a partitioned surface $(\Sigma,P)$.
\begin{itemize}
\item A {\em $P$-separating curve} is a simple closed curve $\gamma$ so that $[\gamma]=0$ in $\HH_1^P(\Sigma)$.
Equivalently, $\gamma$ is a separating curve and for any boundary components $b_1$ and $b_2$ with 
$b_1,b_2 \in p$ for some $p \in P$, the curve $\gamma$ does not separate $b_1$ from $b_2$.
\item A twist about a {\em $P$-bounding pair} equals $T_{\gamma_1} T_{\gamma_2}^{-1}$ for disjoint, nonisotopic simple closed
curves $\gamma_1$ and $\gamma_2$ so that for some choice of orientations we have $[\gamma_1]=[\gamma_2]$ (as elements of $\HH_1^P(\Sigma)$).
\end{itemize}
\end{definition}

\noindent
These are all elements of $\Torelli(\Sigma,P)$:

\begin{lemma}
\label{lemma:separatingbounding}
Fix a partitioned surface $(\Sigma,P)$, and let $f \in \Mod(\Sigma)$ be a twist about either a separating curve or a bounding pair.  Then
$f \in \Torelli(\Sigma,P)$ if and only if $f$ is a twist about a $P$-separating curve or a $P$-bounding pair.
\end{lemma}
\begin{proof}
The $P$-separating curves and $P$-bounding pairs are exactly the separating curves and bounding pairs in $\Sigma$ which remain separating
curves and bounding pairs in any capping of $(\Sigma,P)$.
\end{proof}

\subsection{Morphisms between Torelli groups}
\label{section:torellimorphisms}
In this section, we construct a category $\TSur$ whose objects are partitioned surfaces $(\Sigma,P)$ and whose
morphisms from $(\Sigma_{g_1,n_1},P_1)$ to $(\Sigma_{g_2,n_2},P_2)$ are exactly those embeddings $i \co \Sigma_{g_1,n_1} \hookrightarrow \Sigma_{g_2,n_2}$ 
which induce
morphisms $i_\ast \co \Torelli(\Sigma_{g_1,n_1},P_1) \rightarrow \Torelli(\Sigma_{g_2,n_2},P_2)$.  There is one obvious condition on
such embeddings : for any $P_1$-separating curve $\gamma$, the curve $i(\gamma)$ must be a $P_2$-separating curve.  To translate this into a condition
on the partitions, we need some notation.  Let $S$ be any component of $\Sigma_{g_2,n_2} \setminus \Interior(\Sigma_{g_1,n_1})$.
Observe that $S$ may in fact consist of a single boundary component of $\Sigma_{g_1,n_1}$ which is also a boundary
component of $\Sigma_{g_2,n_2}$.  Let $B_S$ be the boundary components of $\Sigma_{g_2,n_2}$ which lie in
$S$, and let $B_S'$ be the boundary components of $\Sigma_{g_1,n_1}$ which lie in $S$.  Our category is the following:

\begin{definition}
The {\em Torelli surface category} (denoted $\TSur$) is the category whose objects are partitioned surfaces $(\Sigma,P)$ and whose morphisms
from $(\Sigma_{g_1,n_1},P_1)$ to $(\Sigma_{g_2,n_2},P_2)$ are embeddings $\Sigma_{g_1,n_1} \hookrightarrow \Sigma_{g_2,n_2}$
satisfying the following two conditions:
\begin{enumerate}
        \item Each set $B_S'$ is contained in some $p \in P_1$.
        \item Consider $b_1,b_2 \in p$ for some $p \in P_2$.  Assume that $b_1 \in B_S$ for some $S$ and that $b_2 \in B_{S'}$ for some $S' \neq S$.
        Then $B_S' \cup B_{S'}' \subset q$ for some $q \in P_1$.  Less formally, there is a well-defined ``retraction map'' $P_2 \rightarrow P_1$.
\end{enumerate}
\end{definition}

Condition 1 is necessary for all $P_1$-separating curves in
$\Sigma_{g_1,n_1}$ to remain separating curves in $\Sigma_{g_2,n_2}$, and condition 2 is necessary to assure that
they in fact are $P_2$-separating curves.  The following theorem says that these are exactly the morphisms we want: 

\begin{theorem}
\label{theorem:torellimorphisms}
Let $(\Sigma_{g_1,n_1},P_1)$ and $(\Sigma_{g_2,n_2},P_2)$ be partitioned surfaces.  Fix an embedding 
$i \co \Sigma_{g_1,n_1} \hookrightarrow \Sigma_{g_2,n_2}$.  Then
$$i_\ast (\Torelli(\Sigma_{g_1,n_1},P_1)) \subset \Torelli(\Sigma_{g_2,n_2},P_2)$$ 
if and only if $i$ is a morphism of $\TSur$. 
\end{theorem}

\begin{corollary}
\label{corollary:torellifunctor}
$\Torelli$ is a functor from $\TSur$ to the category of groups and homomorphisms.
\end{corollary}

\noindent
We begin by proving a special case.

\begin{proofof}{reverse implication of Theorem \ref{theorem:torellimorphisms} when $\boldsymbol{n_2=0}$}
In this case, condition 2 is vacuous.  By Lemma \ref{lemma:embeddingsarecappings}, there is a partition $P_1'$ of the boundary components
of $\Sigma_{g_1,n_1}$ so that $i \co \Sigma_{g_1,n_1} \hookrightarrow \Sigma_{g_2}$ is a capping of $(\Sigma_{g_1,n_1},P_1')$.  
Observe that condition 1 says that for all
$p' \in P_1'$, there is some $p \in P_1$ with $p' \subset p$.  Defining $H \subset \HH_1^{P_1}(\Sigma)$ to equal
\begin{align*}
\langle \{[h] \text{ $|$ } & \text{$h$ is a simple closed curve or a properly embedded} \\
                           & \text{arc between $P_1'$-adjacent boundary components}\} \rangle,
\end{align*}
we conclude that $H$ surjects onto $\HH_1^{P_1'}(\Sigma)$.  Theorem \ref{theorem:torelliwelldefined} therefore
implies that 
$$\Torelli(\Sigma_{g_1,n_1},P_1) \subset \Torelli(\Sigma_{g_1,n_1},P_1') = i_{\ast}^{-1} (\Torelli(\Sigma_{g_2})),$$
whence the theorem.
\end{proofof}

\noindent
We now prove the general case.

\Figure{figure:torellimorphismproof}{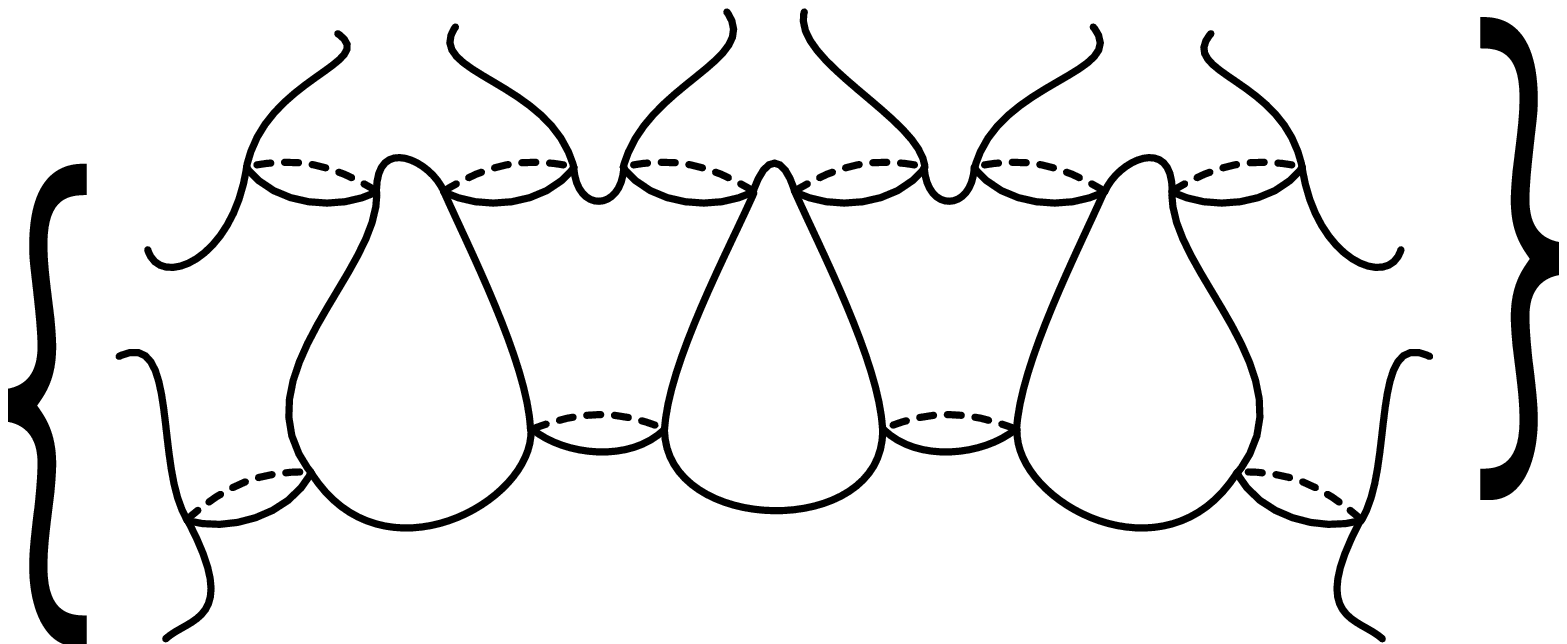}{Portion of a component of $\Sigma_{g_3} \setminus \Interior(\Sigma_{g_1,n_1})$}

\begin{proofof}{\ref{theorem:torellimorphisms} in the general case}
As was mentioned before the theorem, if condition 1 is not satisfied, then there is some $P_1$-separating curve $\gamma$ in
$\Sigma_{g_1,n_1}$ so that $i(\gamma)$ is not a separating curve, and if condition 2 is not satisfied, then there
is some $P_1$-separating curve $\gamma$ in $\Sigma_{g_1,n_1}$ so that $i(\gamma)$ is not a $P_2$-separating curve.  Lemma
\ref{lemma:separatingbounding} therefore implies the necessity of the 2 conditions.  We must prove their
sufficiency.

Let $j \co \Sigma_{g_2,n_2} \hookrightarrow \Sigma_{g_3}$ be a capping of $(\Sigma_{g_2,n_2},P_2)$.  We claim that
$j \circ i \co \Sigma_{g_1,n_1} \hookrightarrow \Sigma_{g_3}$ satisfies the conditions of the theorem (for the partition $P_1$
of the boundary components of $\Sigma_{g_1,n_1}$).  Indeed, if $S$ is
a component of $\Sigma_{g_3} \setminus \Interior(\Sigma_{g_1,n_1})$, then $B_S'$ is a union of $B_{T}'$ for certain
components $T$ of $\Sigma_{g_2,n_2} \setminus \Interior(\Sigma_{g_1,n_1})$.  By condition 1, each $B_{T}'$ is contained
in some $p \in P_1$.  It is not hard to see from the definition
of a capping (see Figure \ref{figure:torellimorphismproof}) that if $T$ and $T'$ are 2 components 
of $\Sigma_{g_2,n_2} \setminus \Interior(\Sigma_{g_1,n_1})$
so that $B_{T}' \cup B_{T'}' \subset B_S'$, then there are components $T=T_1,\ldots,T_k=T'$ of
$\Sigma_{g_2,n_2} \setminus \Interior(\Sigma_{g_1,n_1})$ so that for $1 \leq i < k$ there is some
$b_i \in B_{T_i}$ and $b_i' \in B_{T_{i+1}}$ so that $b_i,b_i' \in q$ for some $q \in P_2$.  By condition 2,
for $i=1,\ldots,k-1$ we have that $B_{T_i}' \cup B_{T_{i+1}}'$ is a subset of a single partition element of $P_1$.  We conclude
that $B_T'$ and $B_{T'}'$ are contained in the same partition element of $P_1$.  This implies that condition 1 holds.  Since condition
2 is vacuous, this implies that $j \circ i$ does indeed satisfy the conditions of the theorem.

The special case of the theorem proven above therefore implies that 
$$j_{\ast} \circ i_{\ast} (\Torelli(\Sigma_{g_1,n_1},P_1)) \subset \Torelli(\Sigma_{g_3}).$$
Since $\Torelli(\Sigma_{g_2,n_2},P_2) = j_{\ast}^{-1} (\Torelli(\Sigma_{g_3}))$, we conclude that 
$$i_\ast (\Torelli(\Sigma_{g_1,n_1},P_1)) \subset \Torelli(\Sigma_{g_2,n_2},P_2),$$
as desired.
\end{proofof}

\subsection{Restriction properties of $\Torelli$}
\label{subsection:restriction}
We now prove the following generalization of Theorem \ref{theorem:torelliwelldefined}:

\begin{theorem}
\label{theorem:torellirestriction}
Let $i \co \Sigma_{g_1,n_1} \hookrightarrow \Sigma_{g_2,n_2}$ be any embedding, and let $P_2$ be any partition
of the boundary components of $\Sigma_{g_2,n_2}$.  There is then some partition $P_1$ of the boundary
components of $\Sigma_{g_1,n_1}$ so that $\Torelli(\Sigma_{g_1,n_1},P_1) = i_{\ast}^{-1}(\Torelli(\Sigma_{g_2,n_2},P_2)$.
\end{theorem}
\begin{proof}
Let
$$j \co \Sigma_{g_2,n_2} \hookrightarrow \Sigma_{g_3}$$
be a capping of $(\Sigma_{g_2,n_2},P_2)$.  By Lemma \ref{lemma:embeddingsarecappings}, there is some partition $P_1$
so that
$$j \circ i \co \Sigma_{g_1,n_1} \hookrightarrow \Sigma_{g_3}$$
is a capping of $(\Sigma_{g_1,n_1},P_1)$.  It is not hard to show that 
$$i \co \Sigma_{g_1,n_1} \hookrightarrow \Sigma_{g_2,n_2}$$
satisfies the conditions of Theorem \ref{theorem:torellimorphisms} (for the partitions $P_i$ of the boundary
components of $\Sigma_{g_i,n_i}$), and hence $i_{\ast}$ induces a morphism from $\Torelli(\Sigma_{g_1,n_1},P_1)$
to $\Torelli(\Sigma_{g_2,n_2},P_2)$.  The theorem is then an immediate corollary
of Theorem \ref{theorem:torelliwelldefined} applied to both $j$ and $i \circ j$.
\end{proof}

\section{A Birman exact sequence for the Torelli groups}
\label{section:torellisequence}
In this section, we prove Theorem \ref{theorem:torelliexactsequenceprelim}.  Fix a partitioned surface $(\Sigma_{g,n},P)$ with
$n \geq 1$ and $(g,n) \neq 1$.
Let $b$ be a boundary component
of $\Sigma_{g,n}$, and let $p \in P$ be the subset containing $b$.  Regarding $\Sigma_{g,n-1}$
as the surface which results from gluing a disc from $b$, let $P'$ be the induced partition
of the boundary components of $\Sigma_{g,n-1}$.  The embedding $i \co \Sigma_{g,n} \rightarrow \Sigma_{g,n-1}$ is
clearly a morphism in $\TSur$ from $(\Sigma_{g,n},P)$ to $(\Sigma_{g,n-1},P')$, and hence
there is an induced map $i_\ast \co \Torelli(\Sigma_{g,n},P) \rightarrow \Torelli(\Sigma_{g,n-1},P')$.  Now,
choosing a point $x$ in the glued-in disc and a lift $\tilde{x}$ of $x$ to the unit tangent
bundle $U\Sigma_{g,n-1}$, we discussed in Section \ref{section:birmanexactsequence} the following
two exact sequences:
\begin{align}
&1 \longrightarrow \pi_1(U\Sigma_{g,n-1},\tilde{x}) \longrightarrow \Mod(\Sigma_{g,n}) \longrightarrow \Mod(\Sigma_{g,n-1}) \longrightarrow 1 \label{seq1},\\
&1 \longrightarrow \Z \longrightarrow \pi_1(U\Sigma_{g,n-1},\tilde{x}) \longrightarrow \pi_1(\Sigma_{g,n-1},x) \longrightarrow 1 \label{seq2}.
\end{align} 
Using exact sequence (\ref{seq1}), we see that $\Ker(i_{\ast}) \subset \pi_1(U\Sigma_{g,n-1},\tilde{x})$.  However,
it is rarely all of $\pi_1(U\Sigma_{g,n-1},\tilde{x})$ : for instance, unless $p=\{b\}$, the twist $T_b$
is not a twist about a $P$-separating curve, so it is not in the kernel.  Theorem \ref{theorem:torelliexactsequenceprelim}
says that $i_{\ast}$ is always surjective and also identifies its kernel.  We will prove a slightly more precise
version of Theorem \ref{theorem:torelliexactsequenceprelim}.  Before stating it, we need the following definition:

\begin{definition}
Assume that a group $\Gamma$ splits as $G_1 \oplus G_2$ and that $\phi \co H_1 \rightarrow G_2$ is a homomorphism,
where $H_1$ is a subgroup of $G_1$.  Then the {\em graph of $\phi$} is the subgroup $\{(x,\phi(x)) \text{ $|$ $x \in H_1$}\}$
of $\Gamma$.
\end{definition}

\noindent
We will prove the following:

\begin{theorem}
\label{theorem:torelliexactsequence}
With the notation as above, we have an exact sequence
$$1 \longrightarrow K \longrightarrow \Torelli(\Sigma_{g,n},P) \longrightarrow \Torelli(\Sigma_{g,n-1},P') \longrightarrow 1$$
with $K \subset \pi_1(U\Sigma_{g,n-1},\tilde{x})$ equal to the following:
\begin{itemize}
\item If $p = \{b\}$, then $K = \pi_1(U\Sigma_{g,n-1},\tilde{x})$.
\item If $p \neq \{b\}$, then $\pi_1(U\Sigma_{g,n-1},\tilde{x})$ splits as $\pi_1(\Sigma_{g,n-1},x) \oplus \Z$ and
$K$ equals the graph of a homomorphism $\phi \co K' \rightarrow \Z$ for the kernel $K' \subset \pi_1(\Sigma_{g,n-1},x)$
of the natural map $\pi_1(\Sigma_{g,n-1},x) \rightarrow \HH_1^{P'}(\Sigma_{g,n-1})$.
\end{itemize}
\end{theorem}

\begin{remark}
The splitting of $\pi_1(U\Sigma_{g,n-1})$ from the second part of Theorem \ref{theorem:torelliexactsequence} is
not natural.  If $\overline{\HH}_1^{P'}(\Sigma_{g,n-1}) = \HH_1(\Sigma_{g,n-1})$ (so $K'$ is the commutator subgroup), then
since exact sequence \ref{seq2} describes a central extension, the embedding of
$K'$ into $\pi_1(U\Sigma_{g,n-1})$ is canonical.  However, this need not be the case.
\end{remark}

\noindent
The proof of the first part of this theorem is easy.

\begin{proofof}{\ref{theorem:torelliexactsequence} when $\boldsymbol{p=\{b\}}$}
Let $i \co \Sigma_{g,n} \hookrightarrow \Sigma_{g,n-1}$ be the embedding.  Observe that since $p=\{b\}$, the map
$i$ induces an isomorphism between $\HH_1^P(\Sigma_{g,n})$ and $\HH_1^{P'}(\Sigma_{g,n-1})$.  This
isomorphism is equivariant with respect to the action of the two mapping class groups, so
Theorem \ref{theorem:torelliwelldefined} implies that $f \in \Torelli(\Sigma_{g,n},P)$ if and only if
$i_{\ast}(f) \in \Torelli(\Sigma_{g,n-1},P')$.  This implies that the standard Birman exact
sequence restricts in the indicated manner, as desired.
\end{proofof}

\noindent
The proof of the other part of Theorem \ref{theorem:torelliexactsequence} is somewhat more complicated.  Our proof
makes use of some ideas of van den Berg \cite[Proposition 2.4.1]{vdB}.

\begin{proofof}{\ref{theorem:torelliexactsequence} when $\boldsymbol{p \neq \{b\}}$}
Observe that since $p \neq \{b\}$, we must have $n-1 \geq 1$.  Thus $\pi_1(\Sigma_{g,n-1},x)$ is a free group and
exact sequence (\ref{seq2}) splits, so $\pi_1(U\Sigma_{g,n-1},\tilde{x}) \cong \pi_1(\Sigma_{g,n-1},x) \oplus \Z$.
Fix such a splitting.  We begin with a criterion for an
element of $\Mod(\Sigma_{g,n})$ to lie in $\Torelli(\Sigma_{g,n},P)$.  Let $h \in \HH_1^P(\Sigma_{g,n})$ be  
the homology class of any arc between $b$ and some other boundary component.

\begin{claim}{1}
Consider $f \in \Mod(\Sigma_{g,n})$.  Assume that $i_{\ast}(f) \in \Torelli(\Sigma_{g,n-1},P')$ and
that $f$ fixes the homology class $h$.  Then $f \in \Torelli(\Sigma_{g,n},P)$.
\end{claim}

\begin{claimproof}
Define $L \subset \HH_1^P(\Sigma_{g,n})$ to equal
\begin{align*}
\langle \{[g] \in \HH_1^P(\Sigma_{g,n}) \text{ $|$ } &\text{$g$ is a simple closed curve or a properly embedded} \\
                                                     &\text{arc between $P'$-adjacent boundary components}\} \rangle.
\end{align*}
It is easy to see that $\HH_1^P(\Sigma_{g,n}) = \langle L, h \rangle$, so it is enough to prove that
$f$ fixes $L$.  Observe that
$$\HH_1^{P'}(\Sigma_{g,n-1}) \cong L / \Span{$b$}.$$
Since $i_{\ast} (f) \in \Torelli(\Sigma_{g,n-1},P')$, it follows that for $g \in L$ we have $f_\ast (g) = g + k \cdot [b]$ for some integer $k$ (the integer $k$ depends on $g$).
However,
$$i_a(h,g) = i_a(f_\ast (h), f_\ast (g)) = i_a(h,g+k \cdot [b]) = i_a(h,g) + k.$$
This implies that $k=0$, as desired.
\end{claimproof}

\noindent
We now examine the manner in which $\pi_1(U\Sigma_{g,n-1},\tilde{x}) \subset \Mod(\Sigma_{g,n})$ acts
on $h$.  Observe that we have an injection
$$j \co \HH_1^{P'}(\Sigma_{g,n-1}) \hookrightarrow \HH_1^P(\Sigma_{g,n}) / \Span{$[b]$}.$$
Since $\Mod(\Sigma_{g,n})$ preserves $[b]$, the action of $\Mod(\Sigma_{g,n})$ on $\HH_1^P(\Sigma_{g,n})$ descends 
to an action on $\HH_1^P(\Sigma_{g,n}) / \Span{$[b]$}$.  Recall that $T_b$ is the generator
of the kernel of exact sequence (\ref{seq2}).  Since $T_b$ 
acts trivially on $\HH_1^P(\Sigma_{g,n}) / \Span{$[b]$}$,
exact sequence (\ref{seq2}) implies that the action of $\pi_1(U\Sigma_{g,n-1},\tilde{x})$ on $\HH_1^P(\Sigma_{g,n}) / \Span{$[b]$}$
descends to an action of $\pi_1(\Sigma_{g,n-1},x)$ on $\HH_1^P(\Sigma_{g,n}) / \Span{$[b]$}$.  Let $\overline{h}$ be the image of $h$
in $\HH_1^P(\Sigma_{g,n}) / \Span{$[b]$}$.  The following claim says that
$\pi_1(\Sigma_{g,n-1},x)$ acts in the most obvious possible way on $\overline{h}$:

\begin{claim}{2}
$\gamma \in \pi_1(\Sigma_{g,n-1},x)$ acts upon $\overline{h}$ by
$$\overline{h} \longmapsto \overline{h} + j([\gamma]).$$
\end{claim}

\begin{claimproof}
Assume first that $\gamma$ is a simple closed curve.  In this case, the lift $\LPtPsh{\gamma} \in \Mod(\Sigma_{g,n})$
(see Figure \ref{figure:birmanexactsequence}.b) equals $T_{\tilde{\gamma}_1} T_{\tilde{\gamma}_2}^{-1}$ for two simple closed curve $\tilde{\gamma}_1$ and
$\tilde{\gamma}_2$ satisfying
$$[\tilde{\gamma}_2] = [\tilde{\gamma}_1] + [b].$$
Now, we know that $i_a(h,[b])=1$.  The claim then follows from an easy calculation, using the fact that Dehn twists act as transvections on homology.

To prove the general case, observe that by the proof of the case $p = \{b\}$, we know that $\pi_1(U\Sigma_{g,n},\tilde{x})$ acts as the identity on
the submodule of $\HH_1^P(\Sigma_{g,n})$ generated by the homology classes of simple closed curves.  Now,
since simple closed curves generate $\pi_1(\Sigma_{g,n-1},x)$, we can write $\gamma = \gamma_1 \gamma_2 \cdots \gamma_k$, where
the $\gamma_i$ are simple closed curves.  The previous case (plus the observation at the beginning of
this paragraph) then implies that
\begin{align*}
\gamma(\overline{h}) &= (\gamma_1 \gamma_2 \cdots \gamma_k)(\overline{h}) = (\gamma_1 \gamma_2 \cdots \gamma_{k-1})(\overline{h} + j([\gamma_k])) \\
                     &= (\gamma_1 \gamma_2 \cdots \gamma_{k-2})(\overline{h} + j([\gamma_{k-1}]) + j([\gamma_{k}])) \\
                     &= \cdots = \overline{h} + j([\gamma_1]) + \cdots + j([\gamma_k]) = \overline{h} + j([\gamma]),
\end{align*}
as desired.
\end{claimproof}

\noindent
We now prove that $K$ equals the graph of some subgroup $K' \subset \pi_1(\Sigma_{g,n-1},x)$.

\begin{claim}{3}
There exists a subgroup $K'$ of $\pi_1(\Sigma_{g,n-1},x)$ and a homomorphism $\phi \co K' \rightarrow \Z$ so
that $K$ is the graph of $\phi$.
\end{claim}

\begin{claimproof}
Let $K'$ equal the projection of $K \subset \pi_1(\Sigma_{g,n-1},x) \oplus \Z$ to $\pi_1(\Sigma_{g,n-1},x)$.  To prove
that $K$ is the graph of a homomorphism $\phi \co K' \rightarrow \Z$, it is enough to show that each element of 
$K'$ is the image of exactly one element of $K$.  In other words, we must prove that $\Z \cap K = 1$.  
This follows from the fact that $T_{b}^m (h) = h+m[b] \neq h$ if $m \neq 0$.
\end{claimproof}

\noindent
We now identify $K'$.

\begin{claim}{4}
$K'$ equals the kernel of the natural map $\pi_1(\Sigma_{g,n-1},x) \rightarrow \HH_1^{P'}(\Sigma_{g,n-1})$.
\end{claim}

\begin{claimproof}
Let $K''$ be the kernel of the natural map $\pi_1(\Sigma_{g,n-1},x) \rightarrow \HH_1^{P'}(\Sigma_{g,n-1})$, and consider 
$f \in \pi_1(U\Sigma_{g,n-1},\tilde{x})$.  Let $\overline{f} \in \pi_1(\Sigma_{g,n-1},x)$ be the projection of $f$.  To
show that $K'=K''$, we need to show that there exists some $m \in \Z$ so that $f T_{b}^m \in \Torelli(\Sigma_{g,n},P)$
if and only if $\overline{f} \in K''$.  By Claim 1, $f T_{b}^m \in \Torelli(\Sigma_{g,n},P)$ if and only if
$f T_{b}^m$ fixes the homology class $h$.  Now, $T_b^m(h) = h + m[b]$.  It follows that we can find an $m$ such
that $f T_{b}^m$ fixes the homology class $h$ if
and only if $f$ fixes $\overline{h} \in \HH_1^P(\Sigma_{g,n}) / \Span{$[b]$}$ (see the discussion before
Claim 2).  By Claim 2, this will be true if and only if $\overline{f} \in K''$, as desired.
\end{claimproof}

\noindent
We finish by proving that the map $\Torelli(\Sigma_{g,n},P) \rightarrow \Torelli(\Sigma_{g,n-1},P')$
is surjective.

\begin{claim}{5}
Let $\overline{f} \in \Torelli(\Sigma_{g,n-1},P')$.  There then exists some
$f \in \Torelli(\Sigma_{g,n},P)$ so that $i_\ast (f) = \overline{f}$.
\end{claim}

\begin{claimproof}
Let $f \in \Mod(\Sigma_{g,n})$ be any lift of $\overline{f}$.  Since the space of all embeddings of
$h \cup b$ into $\Sigma_{g,n-1}$ which fix the endpoint of $h$ not on $b$ is connected, we can
assume that $f$ fixes $h$.  Claim 5 then tells us that $f \in \Torelli(\Sigma_{g,n},P)$, as desired.
\end{claimproof}

\noindent
This completes the proof of the theorem.
\end{proofof}

\subsection{An addendum to Theorem \ref{theorem:torelliexactsequence}}

\Figure{figure:torellisequenceproof}{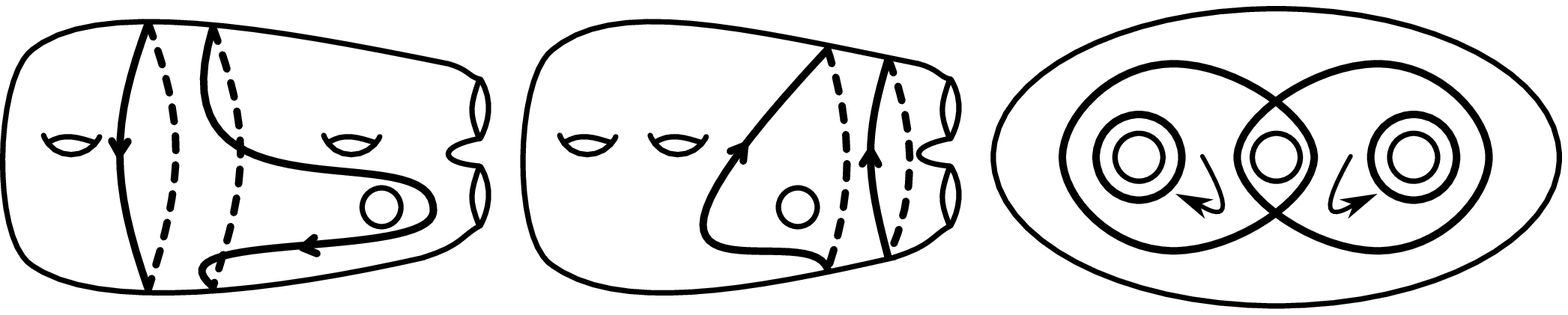}{a. $\LPtPsh{\eta}$ for a simple closed curve 
$\eta$ which cuts off a one-holed torus \CaptionSpace b. $\LPtPsh{\eta}$ for a simple
closed curve $\eta$ which cuts off a set of boundary components \CaptionSpace c. $\LPtPsh{\gamma^1}$ and $\LPtPsh{\gamma^2}$ for
simple closed curves $\gamma^1$ and $\gamma^2$ so that $\gamma^1 \cap \gamma^2 = \{x\}$ and so that a regular neighborhood
of $\gamma^1 \cup \gamma^2$ is homeomorphic to $\Sigma_{0,3}$.}

We now prove the following lemma, whose proof uses some ideas from the proof of \cite[Proposition 2.4.2]{vdB}:

\begin{lemma}
\label{lemma:torellikernel}
Let
$$1 \longrightarrow K \longrightarrow \Torelli(\Sigma_{g,n},P) \longrightarrow \Torelli(\Sigma_{g,n-1},P') \longrightarrow 1$$
be as in Theorem \ref{theorem:torelliexactsequence}.  If $g \geq 1$, then $K$ is in the subgroup
of $\Torelli(\Sigma_{g,n},P)$ generated by twists about $P$-separating curves and twists about $P$-bounding pairs.  If $g = 0$, then
$K$ is in the subgroup generated by twists about $P$-separating curves, twists about $P$-bounding pairs, and commutators
of simply intersecting pairs (see the introduction for the definition of a commutator of a simply intersecting pair).
\end{lemma}
\begin{proof}
As in Theorem \ref{theorem:torelliexactsequence}, let $b$ be the boundary component to which we are attaching a disc and $p \in P$
be the partition element containing $b$.  If $p = \{b\}$, then the lemma is trivial, so we assume that $p \neq \{b\}$.
This implies that $K$ is isomorphic to the kernel of the natural map $\pi_1(\Sigma_{g,n-1}) \rightarrow \HH_1^{P'}(\Sigma_{g,n-1})$.
This kernel is generated by the commutator subgroup $[\pi_1(\Sigma_{g,n-1}),\pi_1(\Sigma_{g,n-1})]$ plus
the set of all simple closed curves which cut off subsets $q \in P'$ of boundary components.  We will
prove that these generators lie in the indicated subgroup of $\Torelli(\Sigma_{g,n},P)$.

We begin with the commutator subgroup.  If $g \geq 1$, then Lemma \ref{lemma:commutatorgenerators}
from the appendix says that the commutator subgroup is generated by $[\gamma^1,\gamma^2]$, where
$\gamma^1,\gamma^2 \in \pi_1(\Sigma_{g,n-1})$ are simple closed curves which only intersect at the
basepoint and where a regular neighborhood of $\gamma^1 \cup \gamma^2$ is homeomorphic to $\Sigma_{1,1}$.
This implies that $[\gamma^1,\gamma^2]$ is homotopic to a simple closed curve $\eta$
which cuts off a one-holed torus.  Replacing $\eta$ by its inverse if necessary, the lift $\LPtPsh{\eta}$ 
then equals $T_{\tilde{\eta}_1} T_{\tilde{\eta}_2}^{-1}$ for simple closed curves
$\tilde{\eta}_1$ and $\tilde{\eta}_2$ like those depicted in Figure \ref{figure:torellisequenceproof}.a.  Observe
that $T_{\tilde{\eta}_1}$ is a twist about a $P$-separating curve and that $T_b T_{\tilde{\eta}_2}^{-1}$ is a
twist about a $P$-bounding pair.  We conclude that the element of $\Torelli(\Sigma_{g,n},P)$ associated
to $[\gamma^1,\gamma^2]$ equals 
$$\LPtPsh{\eta} T_b = T_{\tilde{\eta}_1} (T_b T_{\tilde{\eta}_2}^{-1}),$$
which is in the desired subgroup.

If $g = 0$, then pairs of curves like those in the previous paragraph do not exist.  However, 
it is immediate that the commutator subgroup is generated by $[\gamma^1,\gamma^2]$, where
$\gamma^1,\gamma^2 \in \pi_1(\Sigma_{g,n-1})$ are simple closed curves which only intersect at the 
basepoint and where a regular neighborhood of $\gamma^1 \cup \gamma^2$ is homeomorphic to $\Sigma_{0,3}$.  In this
case, for $1 \leq i \leq 2$ we have that $\LPtPsh{\gamma^i} =  T_{\tilde{\gamma}^i_1}^{\pm 1} T_{\tilde{\gamma}^i_2}^{\mp 1}$
for curves like those depicted in Figure \ref{figure:torellisequenceproof}.c (the signs depend on the orientations on $\gamma^1$ and $\gamma^2$).
Observe that
$$[T_{\tilde{\gamma}^1_1}^{\pm 1} T_{\tilde{\gamma}^1_2}^{\mp 1}, T_{\tilde{\gamma}^2_1}^{\pm 1} T_{\tilde{\gamma}^2_2}^{\mp 1}]
=[T_{\tilde{\gamma}^1_1}^{\pm 1}, T_{\tilde{\gamma}^2_1}^{\pm 1}].$$
Using the commutator identities $[g_1^{-1},g_2] = [g_2,g_1]^{g_1^{-1}}$ and $[g_1,g_2^{-1}] = [g_2,g_1]^{g_2^{-1}}$ if
necessary, we see that this is a commutator of a simply intersecting pair, thus proving that the element of
$\Torelli(\Sigma_{g,n},P)$ associated to $[\gamma^1,\gamma^2]$ lies in the desired subgroup.

We conclude by considering a simple closed curve $\eta \in \pi_1(\Sigma_{g,n-1},x)$ which cuts
off a subset $q \in P'$ of boundary components.  Reversing the orientation of $\eta$ if necessary, 
the associated element $\LPtPsh{\eta}$ of $\Mod(\Sigma_{g,n})$ equals
$T_{\tilde{\eta}_1} T_{\tilde{\eta}_2}^{-1}$ for the curves $\tilde{\eta}_1$ and $\tilde{\eta}_2$ pictured
in Figure \ref{figure:torellisequenceproof}.b.
There are two cases.  The first case is $q \neq p \setminus \{b\}$.  In this case, $T_{\tilde{\eta}_1}$ is a $P$-separating curve
and $T_b T_{\tilde{\eta}_2}^{-1}$ is a $P$-bounding pair.  Hence the element of $\Torelli(\Sigma_{g,n},P)$
associated to $\eta$ equals
$$\LPtPsh{\eta} T_b = T_{\tilde{\eta}_1} (T_b T_{\tilde{\eta}_2}^{-1}),$$
which is in the desired subgroup.  The other case is $q = p \setminus \{b\}$.  In this case,
$T_{\tilde{\eta}_2}$ is a $P$-separating curve
and $T_{\tilde{\eta}_1} T_b^{-1}$ is a $P$-bounding pair.  Hence the element of
$\Torelli(\Sigma_{g,n},P)$ associated to $\eta$ is
$$\LPtPsh{\eta} T_b^{-1} = (T_{\tilde{\eta}_1} T_b^{-1}) T_{\tilde{\eta}_2}^{-1},$$
which again is in the desired subgroup.
\end{proof}

\section{Generators for the Torelli groups}
\label{section:torelligenerators}
In this section, we will prove Theorems \ref{theorem:torelligenerators} and \ref{theorem:torelligenerators2}.  
Our main tool will be Theorem \ref{theorem:armstrong}.  We will apply this theorem to the action
of the Torelli group on the complex of curves, which is a simplicial complex introduced
by Harvey \cite{Hv} to encode the combinatorial topology of a surface.

\begin{definition}
The {\em complex of curves} $\Curves(\Sigma_{g})$ is the simplicial complex whose
simplices are sets $\{\gamma_1,\ldots,\gamma_k\}$ of isotopy classes of nontrivial (that is,
not isotopic to a point) simple closed curves which can be realized disjointly.
\end{definition}

\noindent
The following theorem is due to Harer \cite{Ha}; see \cite{I1} for an alternate
proof:

\begin{theorem}
\label{theorem:curvessimple}
$\Curves(\Sigma_{g})$ is $2g-3$-connected.  In particular, $\Curves(\Sigma_{g})$
is simply-connected for $g \geq 2$.
\end{theorem}

\noindent
Observe that $\Torelli(\Sigma_{g})$ acts on $\Curves(\Sigma_{g})$.  A theorem
of Ivanov \cite[Theorem 1.2]{IvanovBook} implies the following:

\begin{theorem}
\label{theorem:torellinice}
For all $g$, the group $\Torelli(\Sigma_{g})$ acts on $\Curves(\Sigma_{g})$ without rotations.
\end{theorem}

\noindent
The following theorem will be the key to our argument:

\begin{theorem}
\label{theorem:curvestorellisimple}
For $g \geq 2$, the complex $\Curves(\Sigma_{g}) / \Torelli(\Sigma_{g})$ is simply-connected.
\end{theorem}

\noindent
We postpone the proof of Theorem \ref{theorem:curvestorellisimple} until Section \ref{section:curvestorelli}.  Instead,
we use it to prove Theorems \ref{theorem:torelligenerators} and \ref{theorem:torelligenerators2}.

\begin{proofof}{\ref{theorem:torelligenerators} and \ref{theorem:torelligenerators2}}
Observe first that Theorem \ref{theorem:torelligenerators2} follows from repeated applications
of Lemma \ref{lemma:torellikernel}.  We will prove Theorem \ref{theorem:torelligenerators} by induction on $g$.  
The base case is $g=1$.  Using Lemma \ref{lemma:torellikernel}, we can assume that $n=0$, in which case the theorem
is trivial.

Now consider a partitioned surface $(\Sigma_{g,n},P)$ with $g>1$.  By repeated use of Lemma \ref{lemma:torellikernel}, 
we can reduce to the case that $n=0$ (and hence forget about $P$).
By Theorems \ref{theorem:curvessimple}, \ref{theorem:curvestorellisimple}, and
\ref{theorem:torellinice}, we can use Theorem \ref{theorem:armstrong} to conclude that $\Torelli(\Sigma_{g})$ is generated
by the subgroups $\Gamma_{\gamma}$ of $\Torelli(\Sigma_{g})$ stabilizing simple closed curves $\gamma$.  
If $\gamma$ is a separating curve which separates $\Sigma_g$ into
two surfaces $\Sigma_{h_1,1}$ and $\Sigma_{h_2,1}$ with $h_1+h_2=g$, then we have an exact sequence
$$1 \longrightarrow \Span{$T_{\gamma}$} \longrightarrow \Torelli(\Sigma_{h_1,1}) \oplus \Torelli(\Sigma_{h_2,1}) \longrightarrow \Gamma_{\gamma} \longrightarrow 1.$$
By induction, both $\Sigma_{h_1,1}$ and $\Sigma_{h_2,1}$ are generated by twists about separating curves and bounding pairs,
so we conclude that $\Gamma_{\gamma}$ is also generated by such elements, as desired.  If instead $\gamma$ is a nonseparating
curve, let $N$ be a small open regular neighborhood of $\gamma$.  The surface $S=\Sigma_g \setminus N$ is then a genus $g-1$ surface
with 2 boundary components $b_1$ and $b_2$, and additionally the inclusion $S \hookrightarrow \Sigma_g$
is a capping of $(S,P)$ for the partition $P=\{\{b_1,b_2\}\}$ of the boundary components of $S$.  By induction, the group
$\Torelli(S,P)$ is generated by $P$-separating curves and $P$-bounding pairs.  Now, 
we have an exact sequence
$$1 \longrightarrow \Span{$T_{b_1} T_{b_2}^{-1}$} \longrightarrow \Torelli(S,P) \longrightarrow \Gamma_{\gamma} \longrightarrow 1.$$
Since the $P$-separating curves and $P$-bounding pairs in $\Torelli(S,P)$ project to separating curves and bounding pairs in
$\Gamma_{\gamma}$, we conclude that $\Gamma_{\gamma}$ is generated by such elements, and we are done.
\end{proofof} 

\section{The connectivity of $\Curves(\Sigma_{g})/\Torelli(\Sigma_{g})$}
\label{section:curvestorelli}
In this section, we prove Theorem \ref{theorem:curvestorellisimple}.  First, in
Section \ref{section:concrete} we interpret a large subcomplex of $\Curves(\Sigma_g) / \Torelli(\Sigma_g)$
in terms of $\HH_1(\Sigma_g)$.  Next, in Section \ref{section:avoidhomology} we show how to homotope
loops in $\Curves(\Sigma_g) / \Torelli(\Sigma_g)$ so that they ``avoid'' a fixed homology class.  In
Section \ref{section:genus3} this leads quickly to a proof of Theorem \ref{theorem:curvestorellisimple}
when $g \geq 3$.  In the genus $2$ case, the additional argument needed is given
in Section \ref{section:genus2}.

\subsection{A concrete description of a subcomplex of $\Curves(\Sigma_g) / \Torelli(\Sigma_g)$}
\label{section:concrete}
In this section, we give a concrete description of the quotient by $\Torelli(\Sigma_g)$
of the following subcomplex of $\Curves(\Sigma_g)$:

\begin{definition}
The {\em nonseparating complex of curves} on $\Sigma_g$, which we will denote
$\CNoSep(\Sigma_g)$, is the subcomplex of $\Curves(\Sigma_g)$ whose simplices
are sets $\{\gamma_1,\ldots,\gamma_k\}$ of isotopy classes of simple closed curves on
$\Sigma_g$ so that $\Sigma_g \setminus (\gamma_1 \cup \cdots \cup \gamma_k)$ is connected.
\end{definition}

\noindent
The importance of this subcomplex for us comes from the following:
\begin{lemma}
\label{lemma:cnoseptorellibigenough}
The inclusion $\CNoSep(\Sigma_g) / \Torelli(\Sigma_g) \hookrightarrow \Curves(\Sigma_g) / \Torelli(\Sigma_g)$ induces
a surjection on $\pi_1$.
\end{lemma}
\begin{proof}
Fix a base point in $\Curves(\Sigma_g) / \Torelli(\Sigma_g)$
at a vertex corresponding to a nonseparating curve, and consider a based loop $\ell \in \pi_1(\Curves(\Sigma_g) / \Torelli(\Sigma_g))$.
We can assume that $\ell$ is a simplicial path in the $1$-skeleton.  Lift $\ell$ (one edge at a time) to a path $\tilde{\ell}$
in $\Curves(\Sigma_{g})$.  Assume that $\tilde{\ell}$ contains a subpath
of the form $\gamma_1-\gamma_2-\gamma_3$, where $\gamma_2$ is a separating curve.
If $\gamma_1$ and $\gamma_3$ lie on different sides of $\gamma_2$, then they are disjoint and we can replace
$\gamma_1-\gamma_2-\gamma_3$ with $\gamma_1-\gamma_3$.  Otherwise, let $\gamma_2'$ be a nonseparating curve on the side of $\gamma_2$ not
containing $\gamma_1$
and $\gamma_3$.  We can then homotope our segment to $\gamma_1-\gamma_2'-\gamma_3$, eliminating
the separating curve $\gamma_2$.  A similar argument allows us to eliminate any edges corresponding to pairs of
nonseparating curves which together separate the surface.  This allows us to homotope $\tilde{\ell}$ into
the $1$-skeleton of $\CNoSep(\Sigma_g)$.  Projecting this homotopy to $\Curves(\Sigma_g) / \Torelli(\Sigma_g)$,
the based loop $\ell$ is homotoped to a loop coming from a loop in $\CNoSep(\Sigma_g) / \Torelli(\Sigma_g)$, as desired.
\end{proof}

Now, consider a simplex $\{\alpha_1,\ldots,\alpha_k\}$ of $\CNoSep(\Sigma_{g})$.  Observe
that since $\Sigma_g \setminus (\alpha_1 \cup \cdots \cup \alpha_k)$ is connected, we can find 
nonseparating curves $\{\alpha_{k+1},\ldots,\alpha_g,\beta_1,\ldots,\beta_g\}$ so that $i_g(\alpha_i,\alpha_j)=i_g(\beta_i,\beta_j)=0$ and
$i_g(\alpha_i,\beta_j)=\delta_{ij}$ for all $i$ and $j$ and so that $\Sigma_g \setminus (\alpha_1 \cup \beta_1 \cup \cdots \cup \alpha_g \cup \beta_g)$ 
is connected.  If we orient the $\alpha_i$ and the $\beta_j$ correctly, this implies 
that $\{[\alpha_1],[\beta_1],\ldots,[\alpha_g],[\beta_g]\}$ is a symplectic basis for $\HH_1(\Sigma_{g})$.  In particular,
$\{[\alpha_1],\ldots,[\alpha_k]\}$ spans a $k$-dimensional isotropic summand of $\HH_1(\Sigma_{g})$ (a summand $A$
of $\HH_1(\Sigma_g)$ is {\em isotropic} if $i_a(x,y)=0$ for all $x,y \in A$).  Since
we had to put arbitrary orientations on the $\alpha_i$ and the $\beta_j$, this suggests the following definition:

\begin{definition}
Let $V$ be a $\Z$-module with a symplectic form $i(\cdot,\cdot)$.  The {\em complex
of unimodular isotropic lines} in $V$, denoted $\ComplexLines(V)$, is the simplicial complex whose simplices
are sets $\{L_1,\ldots,L_k\}$ of $1$-dimensional summands $L_i$ of $V$ so that
$\Span{$L_1,\ldots,L_k$}$ is a $k$-dimensional isotropic summand of $V$.
\end{definition}

\noindent
Since the Torelli group preserves $\HH_1(\Sigma_{g})$, there is a natural map
$$\pi \co \CNoSep(\Sigma_{g}) / \Torelli(\Sigma_{g}) \longrightarrow \ComplexLines(\HH_1(\Sigma_g)).$$
We will prove the following:

\begin{lemma}
\label{lemma:cnoseptorelliisoset}
$\pi$ is an isomorphism of simplicial complexes.
\end{lemma}
\begin{proof}
We have a series of projections
$$\CNoSep(\Sigma_{g}) \stackrel{\tilde{\pi}}{\longrightarrow} \CNoSep(\Sigma_{g}) / \Torelli(\Sigma_{g}) \stackrel{\pi}{\longrightarrow} \ComplexLines(\HH_1(\Sigma_g)).$$
We must prove that for all simplices $s$ of $\ComplexLines(\HH_1(\Sigma_g))$, there is some simplex $\tilde{s}$
of $\CNoSep(\Sigma_{g})$ so that $\pi \circ \tilde{\pi} (\tilde{s}) = s$, and in addition if $\tilde{s}_1$ and $\tilde{s}_2$ are
simplices of $\CNoSep(\Sigma_{g})$ so that $\pi \circ \tilde{\pi} (\tilde{s}_1) = \pi \circ \tilde{\pi} (\tilde{s}_2)$, then
there is some $f \in \Torelli(\Sigma_{g})$ so that $f(\tilde{s}_1) = \tilde{s}_2$.  We begin with the first assertion.  Let
$s$ be a simplex of $\ComplexLines(\HH_1(\Sigma_g))$.  Pick a symplectic basis $\{a_1,b_1,\ldots,a_g,b_g\}$ so that
$s = \{\Span{$a_1$},\ldots,\Span{$a_h$}\}$.
Lemma \ref{lemma:extendinghomologybasis} from the appendix allows us to realize
this symplectic basis by simple closed curves $\{\alpha_1,\beta_1,\ldots,\alpha_g,\beta_g\}$.
Observe that $\tilde{s}=\{\alpha_1,\ldots,\alpha_h\}$ is
a simplex of $\CNoSep(\Sigma_{g})$ with $\pi \circ \tilde{\pi} (\tilde{s}) = s$.

We now prove the section assertion.  Let $\tilde{s}_1$ and $\tilde{s}_2$ be two simplices of $\CNoSep(\Sigma_{g})$ with 
$\pi \circ \tilde{\pi} (\tilde{s}_1) = \pi \circ \tilde{\pi} (\tilde{s}_2)$.  Let the vertices of the $\tilde{s}_i$ be
$\{\alpha_1^i,\ldots,\alpha_h^i\}$.
Order these and pick orientations so that $[\alpha_j^1] = [\alpha_j^2]$.  Set $a_j = [\alpha_j^1]$, and extend this to a
symplectic basis $\{a_1,b_1\ldots,a_g,b_g\}$ for homology.  Lemma \ref{lemma:extendinghomologybasis} from the appendix allows us
to extend $\{\alpha_1^i,\ldots,\alpha_h^i\}$ to a set of oriented simple closed curves
$\{\alpha_1^i,\beta_1^i,\ldots,\alpha_g^i,\beta_g^i\}$ realizing the homology basis $\{a_1,b_1\ldots,a_g,b_g\}$.
Using the classification of surfaces, there must exist some 
$f \in \Mod(\Sigma_{g})$ so that $f(\alpha_j^1)=\alpha_j^2$ and $f(\beta_j^1)=\beta_j^2$ for all $j$.
Since we have chosen $f$ so that it fixes a basis for homology, it follows that $f \in \Torelli(\Sigma_{g})$.  The
proof concludes with the observation that $f(\tilde{s}_1) = \tilde{s}_2$.
\end{proof}

\noindent
Henceforth we will identify $\CNoSep(\Sigma_g) / \Torelli(\Sigma_g)$ with $\ComplexLines(\HH_1(\Sigma_g))$.

\subsection{Homotoping loops off a homology class}
\label{section:avoidhomology}
In this section, we will show how to homotope (in $\Curves(\Sigma_g) / \Torelli(\Sigma_g)$) curves
which lie in $\ComplexLines(\HH_1(\Sigma_g))$ so that they ``avoid'' a fixed homology class.  To
make this precise, we need the following definition:

\begin{definition}
Let $V$ be a $\Z$-module with a symplectic form $i(\cdot,\cdot)$, and let $W$ be a submodule
of $V$.  We then define $\ComplexLines^W(V)$ to be the full subcomplex of $\ComplexLines(V)$
whose vertices are $1$-dimensional summands $L$ of $V$ so that $L \subset W$.
\end{definition}

\noindent
We will also need the following standard definition from PL-topology:

\begin{definition}
Let $v$ be a vertex of a simplicial complex $X$.  The {\em link} of $v$ in $X$ (denoted
$\Link_v(X)$) is the subcomplex of $X$ whose simplices are simplices $s$ of $X$ so that
$v \notin s$ and so that $s \cup \{v\}$ is a simplex of $X$.
\end{definition}

\noindent
Fix a symplectic basis $\{a_1,b_1,\ldots,a_g,b_g\}$ for $\HH_1(\Sigma_g)$, and
set
$$W = \Span{$a_1,b_1,\ldots,a_{g-1},b_{g-1},a_g$} \subset \HH_1(\Sigma_g).$$
Our lemma is the following:

\begin{lemma}
\label{lemma:avoidhomology}
Let $i \co \ComplexLines(\HH_1(\Sigma_g)) \hookrightarrow \Curves(\Sigma_g) / \Torelli(\Sigma_g)$ be the inclusion.
\begin{itemize}
\item For $g \geq 2$, let $\ell$ be a simplicial arc in $\ComplexLines(\HH_1(\Sigma_g))$ whose endpoints lie
in $\ComplexLines^W(\HH_1(\Sigma_g))$.  There is then a simplicial arc $\ell'$ in $\ComplexLines^W(\HH_1(\Sigma_g))$
so that $i_{\ast}(\ell)$ is homotopic to $i_{\ast}(\ell')$ (fixing the endpoints).
\item For $g \geq 3$ and $L$ any vertex of $\ComplexLines(\HH_1(\Sigma_g))$, let $\ell$ be a simplicial arc
in $\Link_L(\ComplexLines(\HH_1(\Sigma_g)))$ whose endpoints lie in $\ComplexLines^W(\HH_1(\Sigma_g))$.  There
is then a simplicial arc $\ell'$ in $\ComplexLines^W(\HH_1(\Sigma_g)) \cap \Link_L(\ComplexLines(\HH_1(\Sigma_g)))$
so that $i_{\ast}(\ell)$ is homotopic to $i_{\ast}(\ell')$ (fixing the endpoints).
\end{itemize}
\end{lemma}

\noindent
An immediate corollary of the first conclusion of this lemma and Lemma \ref{lemma:cnoseptorellibigenough}
is the following:
\begin{corollary}
\label{corollary:lineswbigenough}
For $g \geq 2$, the inclusion $\ComplexLines^W(\HH_1(\Sigma_g)) \hookrightarrow \Curves(\Sigma_g) / \Torelli(\Sigma_g)$
induces a surjection on $\pi_1$.
\end{corollary}

\begin{remark}
At first glance, it may appear that $\ComplexLines^W(\HH_1(\Sigma_g))$ is contained in the cone of $\Span{$a_g$}$.  However,
while it is true that $i_a(a_g,x)=0$ for all $x \in W$, this does {\it not} imply that $\Span{$x$}$ is adjacent
to $\Span{$a_g$}$ in $\ComplexLines^W(\HH_1(\Sigma_g))$ (for instance, consider $x=2a_1 + a_g$).
\end{remark}

\noindent
Before proving Lemma \ref{lemma:avoidhomology}, we need a definition.

\begin{definition}
Let $L$ be a $1$-dimensional summand of a $\Z$-module $M$ which has a fixed
free basis $\{x_1,\ldots,x_n\}$.  Observe that $L = \Span{$v$}$ for a
primitive vector $v$ (that is, $v$ is not divisible by any integer $n \geq 2$) and
that $v$ is unique up to multiplication by $\pm 1$.  Expand $v$ as
$$v = c_1 x_1 + \cdots + c_n x_n.$$
For any $1 \leq i \leq n$, we then define $\Rk_{x_i}(L)=|c_i|$
(by the observation this is well-defined).  We will call this the {\em $x_i$-rank} of $L$.
\end{definition}

\noindent
The proof of Lemma \ref{lemma:avoidhomology} then goes as follows:

\begin{proofof}{\ref{lemma:avoidhomology}}
The proofs of both parts of the lemma are similar; we will prove the (slightly more difficult) second
part and leave the first to the reader.

The proof will be by induction on
$$R = \Max\{\Rk_{b_g}(A) \text{ $|$ $A$ is a vertex of $\ell$}\}.$$
The case $R=0$ being trivial, we assume that $R > 0$.  We first claim that we can
assume that $\ell$ does not contain two adjacent vertices $A$ and $B$ so that $\Rk_{b_g}(A) = \Rk_{b_g}(B) = R$.
Indeed, assume it does and that $A=\Span{$v_1$}$ and $B = \Span{$v_2$}$.  Trivially we can assume
that $A \neq B$.  Replacing $v_2$ by $-v_2$ if
necessary, it follows that $\Rk_{b_g}(\Span{$v_1 - v_2$}) = 0$.  Figure \ref{figure:genus2}.a then shows that we can homotope
(in $\Curves(\Sigma_g) / \Torelli(\Sigma_g)$) the segment $A-B$ to $A-\Span{$v_1-v_2$}-B$, rendering $A$ and $B$ nonadjacent.

Now consider some segment $X-A-Y$ of $\ell$ so that $\Rk_{b_g}(A)=R$ and $\Rk_{b_g}(X),\Rk_{b_g}(Y) < R$.  Lift
$X$, $A$, $Y$, and $L$ to curves $\tilde{X}, \tilde{A}, \tilde{Y}, \tilde{L} \in \CNoSep(\Sigma_g)$ so that
$\{\tilde{X},\tilde{A},\Tilde{L}\}$ and $\{\tilde{Y},\tilde{A},\tilde{L}\}$ are simplices of $\CNoSep(\Sigma_g)$.  Cutting
$\Sigma_g$ along $\tilde{A} \cup \tilde{L}$, we get a copy of $\Sigma_{g-2,4}$ in which both $\tilde{X}$ and $\tilde{Y}$ are
nonseparating curves.  Using Lemma \ref{lemma:connectednesslemma} from the Appendix, we can find a sequence of nonseparating curves
$$\tilde{X} = \tilde{Z}_1,\tilde{Z}_2,\ldots,\tilde{Z}_k = \tilde{Y}$$
in $\Sigma_{g-2,4}$ so that for $1 \leq i < k$ we have $i_g(\tilde{Z}_i,\tilde{Z}_{i+1}) = 1$.  Gluing the boundary components of
our $\Sigma_{g-2,4}$ back together to recover our original $\Sigma_{g}$, we have obtained a sequence of curves
$$\tilde{X} = \tilde{Z}_1,\tilde{Z}_2,\ldots,\tilde{Z}_k = \tilde{Y}$$
so that for $1 \leq i \leq k$ the set $\{\tilde{Z}_i, \tilde{A}, \tilde{L}\}$ is a simplex of $\CNoSep(\Sigma_g)$ and
so that for $1 \leq i < k$ we have $i_g(\tilde{Z}_i,\tilde{Z}_{i+1}) = 1$.

Let $a$ be a primitive vector in $\HH_1(\Sigma_g)$ so that $A = \Span{$a$}$, and for $1 \leq i \leq k$ let $z_i$ be
a primitive vector in $\HH_1(\Sigma_g)$ so that $Z_i:=\Span{$z_i$}$ is the projection of $\tilde{Z}_i$ to $\ComplexLines(\HH_1(\Sigma_g))$.
By the division algorithm, we can find integers $q_i$ (with $q_1 = q_k = 0$) so that
$$\Rk_{b_g}(\Span{$z_i + q_i a$}) < R.$$
Set $z_i' = z_i + q_i a$ and $Z_i' = \Span{$z_i'$}$ (hence $Z_1'=Z_1=X$ and $Z_k'=Z_k=Y$).  Observe that
for all $1 \leq i \leq k$, the set $\{Z_i', A, L\}$ is still a simplex of
$\ComplexLines(\HH_1(\Sigma_g))$, and in addition we still have for all $1 \leq i < k$ that $i_a(z_i',z_{i+1}')=1$.  We
can therefore lift the $Z_i'$ to nonseparating curves $\tilde{Z}_i'$ in $\Sigma_{g}$ which are disjoint from $\tilde{A}$
and $\tilde{L}$ and which satisfy $i_g(\tilde{Z}_i',\tilde{Z}_{i+1}')=1$ for all $1 \leq i < k$.

For $1 \leq i < k$, let $\tilde{S}_i$ be the boundary component of a regular neighborhood of $\tilde{Z}_i \cup \tilde{Z}_{i+1}$.  Hence
$\tilde{S}_i$ is a separating curve disjoint from $\tilde{Z}_i$, $\tilde{Z}_{i+1}$, $\tilde{A}$, and $\tilde{L}$.
For $1 \leq i < k$, cutting $\Sigma_g$ along $\tilde{S}_i$ and $L$ yields a copy of $\Sigma_{g-2,3}$.  A simple dimension count
shows that there must be some simple closed nonseparating curve $\tilde{B}_i$ in this copy of
$\Sigma_{g-2,3}$ whose homology class lies in $W$ (and, in
particular, the span of whose homology class has $b_g$-rank $0$).  Observe that the path
$$\tilde{Z}_1 - \tilde{A} - \tilde{Z}_k$$
in $\Curves(\Sigma_g)$ is homotopic (fixing the endpoints) to
$$\tilde{Z}_1 - \tilde{S}_1 - \tilde{Z}_2 - \tilde{S}_2 - \cdots - \tilde{S}_{k-1} - \tilde{Z}_k,$$
which is then homotopic to
$$\tilde{Z}_1 - \tilde{B}_1 - \tilde{Z}_2 - \tilde{B}_2 - \cdots - \tilde{B}_{k-1} - \tilde{Z}_k.$$
Projecting this homotopy down to $\Link_L(\Curves(\Sigma_g) / \Torelli(\Sigma_g))$ allows us to homotope
the segment $X - A - Y$ to a new segment which lies in $\Link_L(\ComplexLines(\HH_1(\Sigma_g)))$ and which does
not contain any vertices whose $b_g$-rank is greater than or equal to $R$.  Repeating this process
allows us to remove all vertices of $\ell$ whose $b_g$-rank equals $R$, and we are done
by induction.
\end{proofof}

\subsection{Completing the proof for $g \geq 3$}
\label{section:genus3}
In this section, we complete the proof of Theorem \ref{theorem:curvestorellisimple} for $g \geq 3$.  Define
$$W' = \Span{$a_1,b_1,\ldots,a_{g-1},b_{g-1}$}.$$
We will prove the following:
\begin{lemma}
\label{lemma:lineswpbigenough}
The inclusion $\ComplexLines^{W'}(\HH_1(\Sigma_g)) \hookrightarrow \Curves(\Sigma_g) / \Torelli(\Sigma_g)$ induces
a surjection on $\pi_1$.
\end{lemma}
First, though, we will use Lemma \ref{lemma:lineswpbigenough} to prove the desired
special case of Theorem \ref{theorem:curvestorellisimple}.

\begin{proofof}{\ref{theorem:curvestorellisimple} for $\boldsymbol{g \geq 3}$}
By Lemma \ref{lemma:lineswpbigenough}, it is enough to show that the inclusion
$\ComplexLines^{W'}(\HH_1(\Sigma_g)) \hookrightarrow \Curves(\Sigma_g) / \Torelli(\Sigma_g)$
induces the zero map on $\pi_1$.  Indeed, the inclusion map
$\ComplexLines^{W'}(\HH_1(\Sigma_g)) \hookrightarrow \ComplexLines^{W}(\HH_1(\Sigma_g))$
induces the zero map on $\pi_1$, as $\ComplexLines^{W'}(\HH_1(\Sigma_g))$ can be contracted
in $\ComplexLines^{W}(\HH_1(\Sigma_g))$ to $\Span{$a_g$}$.
\end{proofof}

\noindent
We now prove Lemma \ref{lemma:lineswpbigenough}.

\begin{proofof}{\ref{lemma:lineswpbigenough}}
Fix a basepoint in $\Curves(\Sigma_g) / \Torelli(\Sigma_g)$ which lies in the image
of $\ComplexLines^{W'}(\HH_1(\Sigma_g))$, and consider a based loop $\ell$, which
we can assume to be a simplicial loop in the $1$-skeleton.  By Lemma
\ref{lemma:avoidhomology}, we can assume that $\ell$ lies in the image
of $\ComplexLines^{W}(\HH_1(\Sigma_g))$.  We will prove that $\ell$ can
be homotoped into the image of $\ComplexLines^{W'}(\HH_1(\Sigma_g))$ by
induction on
$$R = \Max\{\Rk_{a_g}(A) \text{ $|$ $A$ is a vertex of $\ell$ }\}.$$
The case $R=0$ being trivial, we assume that $R>0$.  By an argument like that
in the proof of Lemma \ref{lemma:avoidhomology}, we can assume that no two
adjacent vertices in $\ell$ have $a_g$-rank equal to $R$.  Now consider
a subpath $X - A - Y$ of $\ell$ so that $\Rk_{a_g}(A) = R$ and
$\Rk_{a_g}(X),\Rk_{a_g}(Y) < R$.  Since $\Link_A(\ComplexLines(\HH_1(\Sigma_g))$ is
connected (this is our main use of the assumption $g \geq 3$), we can find a path
$$X=Z_1 - Z_2 - \cdots - Z_k = Y$$
in $\Link_A (\ComplexLines(\HH_1(\Sigma_g)))$, which by Lemma \ref{lemma:avoidhomology} we can assume
lies in $\ComplexLines^W(\HH_1(\Sigma_g))$.  Let $a$ be a primitive
vector in $\HH_1(\Sigma_g)$ so that $A = \Span{$a$}$, and for $1 \leq i \leq k$ let $z_i$ be a primitive vector
in $\HH_1(\Sigma_g)$ so that $Z_i = \Span{$z_i$}$.  By the division algorithm, for $1 \leq i \leq k$ there
exists an integer $q_i$ (with $q_1=q_k=0$) so that $\Rk_{a_g}(\Span{$z_i + q_i a$}) < R$.  Setting $Z_i' = \Span{$z_i + q_i a$}$,
the path
$$X = Z_1' - Z_2' - \cdots Z_k' = Y$$
still lies in $\ComplexLines^W(\HH_1(\Sigma_g)) \cap \Link_A (\ComplexLines(\HH_1(\Sigma_g)))$.  We can therefore
homotope $X - A - Y$ to this path, eliminating $A$ without introducing any new vertices whose $a_g$-rank
is greater than or equal to $R$.  Repeating this process, we can eliminate all vertices whose $a_g$-rank equals
$R$, and we are done by induction.
\end{proofof}

\subsection{Completing the proof for $g = 2$}
\label{section:genus2}
In genus $2$, the above proof fails.  To complete the proof
in this case, we introduce one final object.

\begin{definition}
Let $V$ be a $\Z$-module with a symplectic inner product $i(\cdot,\cdot)$.  The
{\em augmented complex of unimodular isotropic lines} in $V$, which we will denote
$\AugComplexLines(V)$, is equal to $\ComplexLines(V)$ with 2-cells attached to
all triangles of the form
$$\Span{$v_1$} - \Span{$\pm v_1 \pm v_2$} - \Span{$v_2$}$$
for all edges $\{\Span{$v_1$}, \Span{$v_2$}\}$ of $\ComplexLines(V)$.  For
any submodule $W$ of $V$, we also define $\AugComplexLines^W(V)$ to be
the full subcomplex of $\AugComplexLines(V)$ spanned by vertices $L$ with
$L \subset W$.
\end{definition}

\Figure{figure:genus2}{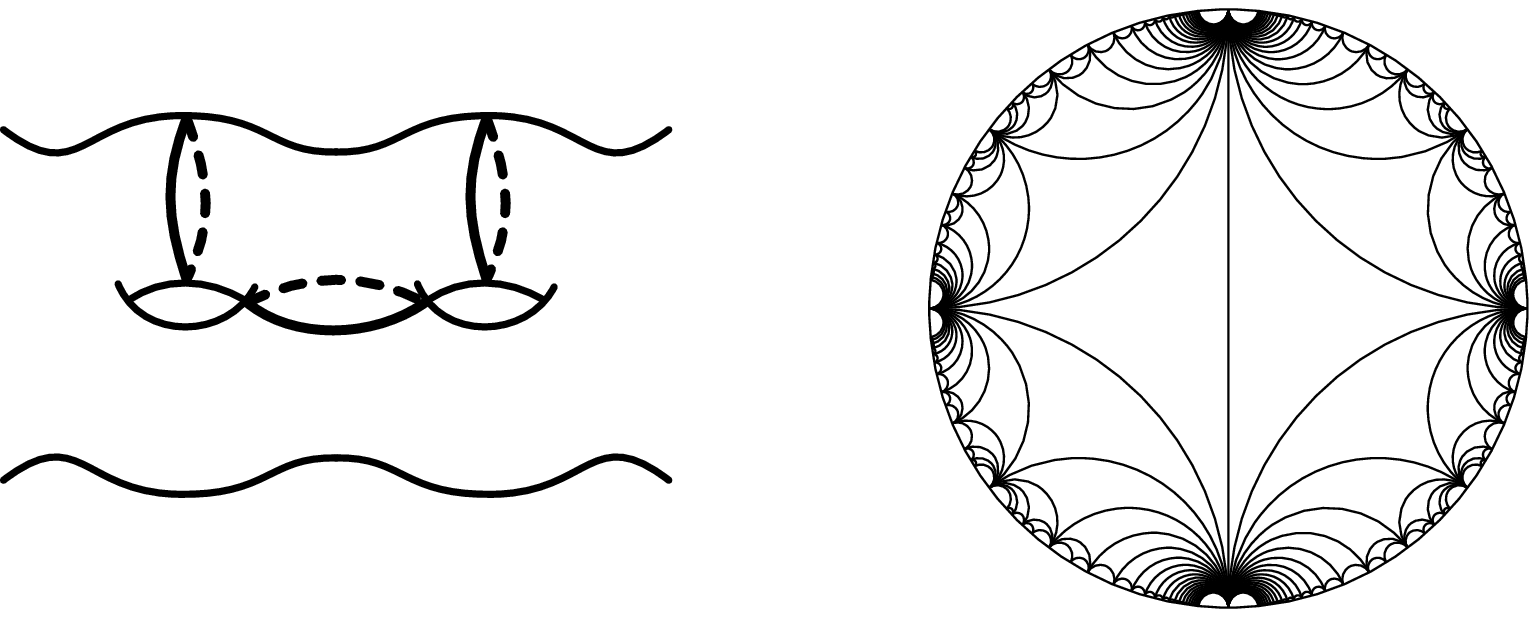}{a. Illustration that triangles $\Span{$v_1$} - \Span{$\pm v_1 \pm v_2$} - \Span{$v_2$}$
in $\ComplexLines(\HH_1(\Sigma_g))$ come from triangles in $\Curves(\Sigma_g)$ \CaptionSpace b. The Farey tessellation}

\noindent
Figure \ref{figure:genus2}.a yields the following lemma:

\begin{lemma}
\label{lemma:augcomplexlinesextends}
The inclusion $\ComplexLines(\HH_1(\Sigma_g)) \hookrightarrow \Curves(\Sigma_g) / \Torelli(\Sigma_g)$
extends to a map $\AugComplexLines(\HH_1(\Sigma_g)) \hookrightarrow \Curves(\Sigma_g) / \Torelli(\Sigma_g)$.
\end{lemma}

\noindent
We now prove the following (recall that the Farey tessellation
of $\HBolic^2$ is the 2-dimensional flag complex whose vertices are elements of
$\Q \cup \{\infty\}$ and where $b/a$ and $d/c$ are adjacent if $|ad-bc|=1$; c.f.\ Figure \ref{figure:genus2}.b):
\begin{lemma}
\label{lemma:augcomplexfarey}
Let $W$ be a maximal isotropic subspace of $\HH_1(\Sigma_2)$.  Then the simplicial complex
$\AugComplexLines^W(\HH_1(\Sigma_2))$ is homeomorphic to the Farey tessellation of $\HBolic^2$ (with
the weak topology).  In particular, $\AugComplexLines^W(\HH_1(\Sigma_2))$ is contractible.
\end{lemma}
\begin{proof}
Identifying $W$ with $\Z^2$, observe that $\AugComplexLines^W(\HH_1(\Sigma_2))$ equals
the following simplicial complex:
\begin{itemize}
\item The vertices are the 1-dimensional summands of $\Z^2$.  These are classified
by their slope in $\Q \cup \{\infty\}$.
\item Two vertices $b/a,d/c \in \Q \cup \infty$ form a 1-simplex if
$\Span{$(a,b),(c,d)$} = \Z^2$.  This is true if and only if the matrix
whose columns are the vectors $(a,b)$ and $(c,d)$ is invertible over $\Z$; i.e.
if and only if the determinant $ad-bc$ equals $\pm 1$.
\item Three vertices $b/a$, $d/c$, and $f/e$ form a 2-simplex if $b/a$ and
$d/c$ form an edge and $(e,f) = \pm (a,b) \pm (c,d)$.  It is easy to see
that this is true if and only if each pair of vertices form an edge.
\end{itemize}
The claim follows.
\end{proof}

\noindent
We now prove Theorem \ref{theorem:curvestorellisimple} when $g=2$.

\begin{proofof}{\ref{theorem:curvestorellisimple} for $\boldsymbol{g=2}$}
Consider a loop $\ell$ in $\Curves(\Sigma_2) / \Torelli(\Sigma_2)$ based at
$\Span{$a_1$} \in \AugComplexLines(\HH_1(\Sigma_2)) \subset \Curves(\Sigma_2) / \Torelli(\Sigma_2)$.  We
can assume that $\ell$ is a simplicial loop in the $1$-skeleton.  By
Lemma \ref{lemma:avoidhomology}, we can also assume that $\ell$ lies in
$\AugComplexLines^{\Span{$a_1,b_1,a_2$}}(\HH_1(\Sigma_2))$.  We will prove that $\ell$ can
be homotoped to a point by induction on the length of $\ell$.  Observe first
that if $\ell$ lies in $\AugComplexLines^{\Span{$a_1,a_2$}}(\HH_1(\Sigma_2))$, then
Lemma \ref{lemma:augcomplexfarey} implies that $\ell$ can be contracted to a point.  Assume,
therefore, that $\ell$ contains a vertex whose $b_1$-rank is nonzero.  Let $\ell'$ be a
subpath of $\ell$ all of whose vertices have $b_1$-rank nonzero and which is maximal
with respect to this property.  In addition, let $x$ and $y$ be the vertices of $\ell$
immediately preceding and immediately succeeding $\ell$'; in other words, $\ell$ contains
the subpath
$$x - \ell' - y.$$
Observe that by the maximality of $\ell'$, both $x$ and $y$ have $b_1$-rank equal to 0; i.e.
$x,y \in \AugComplexLines^{\Span{$a_1,a_2$}}$.  Additionally, since every vertex of $\ell'$
has positive $b_1$-rank, both $x$ and $y$ must have $a_1$-rank $0$.  We conclude that
$x=y=\Span{$a_2$}$.  In other words,
$$x - \ell' - y$$
is a loop.  If $\ell'$ consists of a single vertex, then we can contract this loop and thus
shorten $\ell$, and we are done.  Otherwise, since every vertex of $\ell'$ has $b_1$-rank
greater than $0$, it follows that every vertex of $\ell'$ has $a_1$-rank equal to $0$.  We
conclude that the loop
$$x - \ell' - y$$
lies in $\AugComplexLines^{\Span{$a_2,b_1$}}(\HH_1(\Sigma_2))$.  By Lemma \ref{lemma:augcomplexfarey},
we can contract this loop to a point, thus shortening $\ell$, and we are done.
\end{proofof}

\appendix
\section{Appendix on surface topology}
In this appendix, we will prove three lemmas about the topology of surfaces for which
we are unable to provide appropriate references, though they are certainly known to
the experts.

\subsection{Generators for the commutator subgroup of a surface group}
In this section, we will prove a lemma which provides generators for the commutator subgroup of
a surface group.  It is a generalization of a lemma which appears in a paper of Johnson
\cite[Lemma 7]{J3}.  It also is implicit in the unpublished thesis of van den Berg
\cite[Proposition 2.4.2]{vdB}, though her proof is not quite complete.

\Figure{figure:surfacegroup}{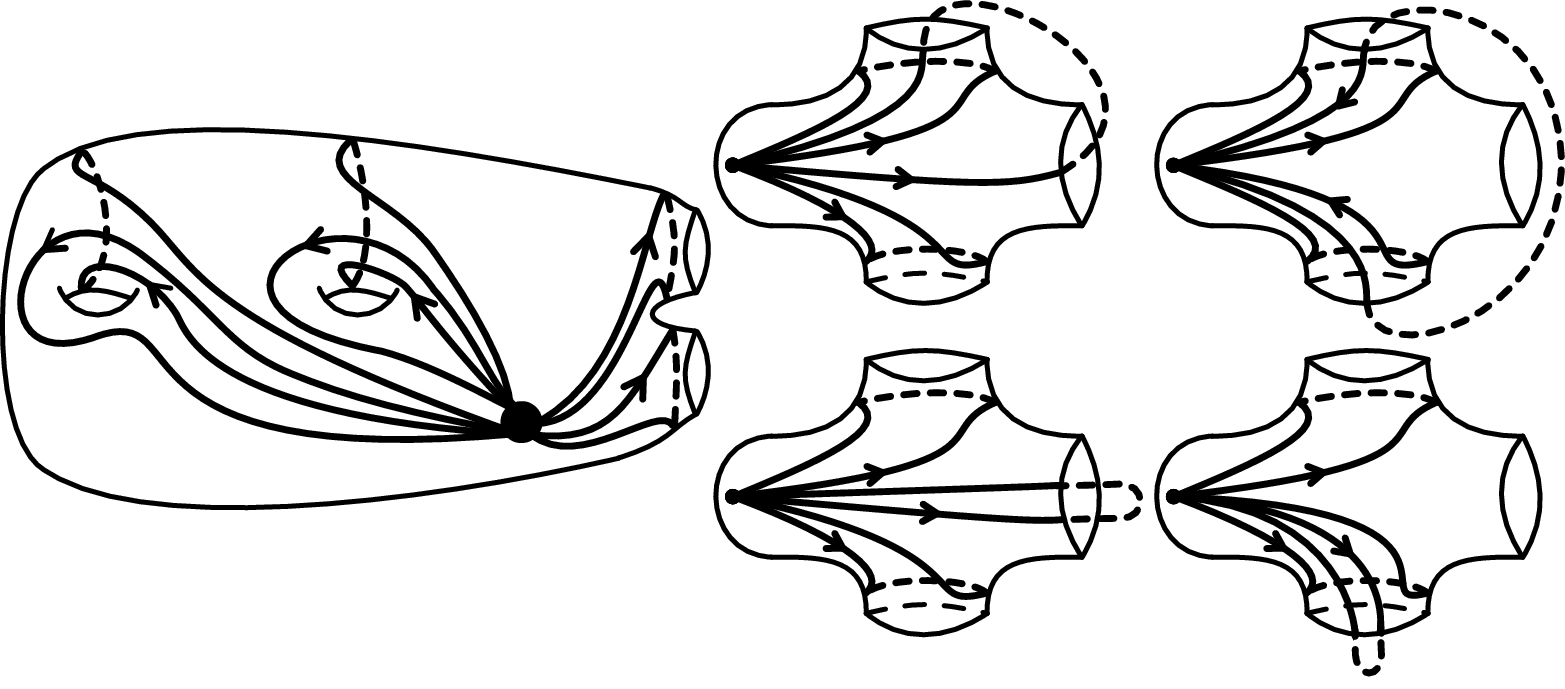}{a. Generators for $\pi_1(\Sigma_{2,2},x)$ \CaptionSpace b. The
four configurations of curves on three-holed spheres needed}

\begin{lemma}
\label{lemma:commutatorgenerators}
Let $g \geq 1$ and $\pi'=[\pi_1(\Sigma_{g,n},x),\pi_1(\Sigma_{g,n},x)]$.  Assume that the basepoint $x$
is in the interior of $\Sigma_{g,n}$.  Then $\pi'$ is
generated by $[\gamma_1,\gamma_2]$, where $\gamma_1,\gamma_2 \in \pi_1(\Sigma_{g,n},x)$ are simple closed curves
so that $\gamma_1 \cap \gamma_2 = \{x\}$ and so that a regular neighborhood of $\gamma_1 \cup \gamma_2$ is homeomorphic to a
one-holed torus.
\end{lemma}
\begin{proof} 
Let $\Gamma$ be the subgroup generated by the indicated elements.  We will first prove
that $\Gamma$ contains $[\gamma_1,\gamma_2]$ for simple closed curves $\gamma_1$ and
$\gamma_2$ so that $\gamma_1 \cap \gamma_2 = \{x\}$ and so that a regular neighborhood
of $\gamma_1 \cup \gamma_2$ is homeomorphic to a three-holed sphere.  There are two
cases.  In the first, one of the $\gamma_i$ (say $\gamma_2$) is nonseparating.  We can
then find a simple closed curve $\gamma_3$ so that $\gamma_3 \cap \gamma_1 = \gamma_3 \cap \gamma_2 = \{x\}$
and so that regular neighborhoods of both $\gamma_1 \gamma_3 \cup \gamma_2$ and $\gamma_3 \cup \gamma_2$ are homeomorphic 
to one-holed tori (see the top portion of Figure \ref{figure:surfacegroup}.b).  We then
have the identity
$$[\gamma_1,\gamma_2] = (\gamma_3 [\gamma_1 \gamma_3,\gamma_2] \gamma_3^{-1}) (\gamma_3 [\gamma_2,\gamma_3] \gamma_3^{-1}).$$
Observe that conjugation does not change the ``type'' of a commutator (this follows
from the identification of $\pi_1(\Sigma_{g,n},x)$ with the kernel in the Birman exact sequence, and is
the reason we require the basepoint to be in the interior).  Hence
the right hand side is in $\Gamma$, so the left hand side is as well.

In the second case, both of the $\gamma_i$ are separating curves.  Reordering the $\gamma_i$ if necessary, we can then
find a simple closed nonseparating curve $\gamma_3$ so that $\gamma_3 \cap \gamma_1 = \gamma_3 \cap \gamma_2 = \{x\}$ and
so that $\gamma_1 \gamma_3$ is nonseparating (see the bottom portion of Figure \ref{figure:surfacegroup}.b; this is where
we use the assumption that $g \geq 1$).  We
again have the identity
$$[\gamma_1,\gamma_2] = (\gamma_3 [\gamma_1 \gamma_3,\gamma_2] \gamma_3^{-1}) (\gamma_3 [\gamma_2,\gamma_3] \gamma_3^{-1}).$$
By the previous case, the right hand side is in $\Gamma$, so the left hand side is as well.

Now let 
$$S=\{\alpha_1,\beta_1,\ldots,\alpha_g,\beta_g,\delta_1,\ldots,\delta_n\}$$ 
be a standard basis for $\pi_1(\Sigma_{g,n},x)$ (see Figure \ref{figure:surfacegroup}.a).  Thus
$\pi'$ is normally generated by $[\gamma_1,\gamma_2]$ for $\gamma_1,\gamma_2 \in S$.  Since we have proven that
every such commutator is in $\Gamma$, we conclude that $\Gamma = \pi'$, as desired.
\end{proof}

\subsection{A connectedness lemma}
In this section, we prove a lemma which allows us to move between simple
closed nonseparating curves in a simple manner.

\begin{lemma}
\label{lemma:connectednesslemma}
For $g \geq 1$ and $n \geq 0$, let $\gamma$ and $\gamma'$ be two simple closed
nonseparating curves in $\Sigma_{g,n}$.  We can then find a sequence $\gamma_1,\gamma_2,\ldots,\gamma_k$
of simple closed nonseparating curves in $\Sigma_{g,n}$ so that $\gamma_1 = \gamma$, $\gamma_k = \gamma'$, and so that for
$1 \leq i <k$ we have $i_g(\gamma_i,\gamma_{i+1})=1$.
\end{lemma}

\begin{remark}
For $n=0$, this lemma is well-known.  However, we need it for
surfaces with boundary, so we include a proof.
\end{remark}

\begin{proofof}{\ref{lemma:connectednesslemma}}
It is well-known that there exists a set $S$ of simple closed curves in $\Sigma_{g,n}$ so
that $\{T_{\delta} \text{ $|$ $\delta \in S$}\}$ generates $\Mod(\Sigma_{g,n})$ and
so that for all $\delta \in S$ we have $i_g(\delta, \gamma) \leq 1$ (for example,
$S$ could be the curves the twists about which form the generating set in \cite{Gervais} and $\gamma$
could be the ``central'' curve $b$ from that paper).  Observe that for all $\delta \in S$, either $T_{\delta}^{\pm 1}(\gamma)=\gamma$
or $i_g(T_{\delta}^{\pm 1}(\gamma), \gamma)=1$.  Now, $\Mod(\Sigma_{g,n})$ acts transitively on
the set of simple closed nonseparating curves, so there exists a sequence of curves $\delta_1,\ldots,\delta_k \in S$
and a sequence of numbers $e_1,\ldots,e_k \in \{\pm 1\}$ so that $T_{\delta_1}^{e_1} \cdots T_{\delta_k}^{e_k} (\gamma) = \gamma'$.
We conclude that after eliminating repetitions, the following is the desired sequence:
$$\gamma, T_{\delta_1}^{e_1}(\gamma), T_{\delta_1}^{e_1} T_{\delta_2}^{e_2}(\gamma), \ldots, T_{\delta_1}^{e_1} T_{\delta_2}^{e_2} \cdots T_{\delta_k}^{e_k}(\gamma).$$
\end{proofof}

\subsection{Realizing homology bases}
In this section, we prove a lemma which allows us to realize symplectic bases for $\HH_1(\Sigma_{g})$ in
a nice manner.

\Figure{figure:extendinghomologybasis}{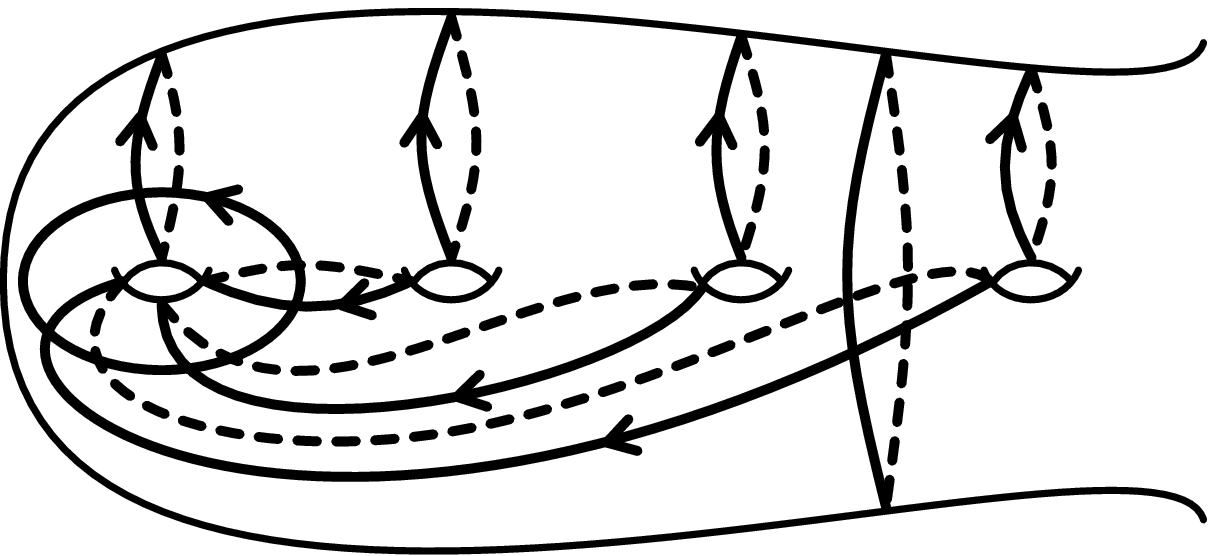}{The various curves needed in the proof of Lemma \ref{lemma:extendinghomologybasis}}

\begin{lemma}
\label{lemma:extendinghomologybasis}
Let $\{a_1,\ldots,a_g,b_1,\ldots,b_g\}$ be a symplectic basis for $\HH_1(\Sigma_{g})$, and for integers
$0 \leq h,k \leq g$ let $\{\alpha_1,\ldots,\alpha_h,\beta_1,\ldots,\beta_k\}$ be a set of oriented
simple closed curves in $\Sigma_{g}$ so that the following hold for all $i$ and $j$ for which
the expressions are defined:
\begin{align*}
[\alpha_i]=a_i \text{ and } [\beta_i]=b_i, \\
i_g(\alpha_i,\alpha_j)=0 \text{ and } i_g(\beta_i,\beta_j)=0, \\
i_g(\alpha_i,\beta_j)=\delta_{ij}.
\end{align*}
Then we can find simple closed curves $\{\alpha_{h+1},\ldots,\alpha_g,\beta_{k+1},\ldots,\beta_g\}$ so that
these expressions continue to hold.
\end{lemma}

\begin{proof}
The proof will be by induction on $g$.  The case $g=0$ is trivial.  Now assume that $g>0$.  If both
$\alpha_1$ and $\beta_1$ are given to us, let $N$ be a small regular neighborhood of
$\alpha_1 \cup \beta_1$.  Observe that $N$ is a copy of $\Sigma_{1,1}$ which is disjoint from the remaining
$\alpha_i$ and $\beta_j$.  Hence $\Sigma_{g} \setminus N$ is a copy of $\Sigma_{g-1,1}$
containing $\{\alpha_2,\ldots,\alpha_h,\beta_2,\ldots,\beta_k\}$.  Let $\Sigma_{g} \setminus N \longrightarrow \Sigma_{g-1}$
be the embedding induced by gluing a disk to the boundary component.  By induction, we can solve
the resulting problem on $\Sigma_{g-1}$, and it is clear that any lift of the solution
to $\Sigma_{g}$ solves the problem there as well.  

We can therefore assume without loss of generality that we are given no $\beta_i$'s.  If no curves are given, 
then it is trivial to find $\alpha_1$, so we can also assume without loss of generality that we are given 
$\{\alpha_1,\ldots,\alpha_h\}$ for
some $h \geq 1$.  Let $\sigma$ be any separating curve dividing $\Sigma_{g}$ into two
subsurfaces $S_1$ and $S_2$ with $\alpha_i \in S_1$ for all $i$.  Arrange the $\alpha_i$ and $\sigma$ in
the pattern indicated in Figure \ref{figure:extendinghomologybasis}, and let $\gamma_2,\ldots,\gamma_h$ and
$\beta'$ be the curves indicated there.  It is clear that with the indicated orientations we have
$$[\gamma_i] = [\alpha_i] - [\alpha_1].$$
Also, $\sigma$ induces a symplectic splitting
$$\HH_1(\Sigma_{g}) = \HH_1(S_1) \oplus \HH_1(S_2).$$
Pick $d \in \Z$ and an irreducible vector $b'' \in \HH_1(S_2)$ so that the projection of $b_1$ to $\HH_1(S_2)$ 
equals $d b''$ (if this projection is $0$, then $d=0$ and $b''$ is arbitrary).  Let $\beta''$ be any simple closed curve in $S_2$ realizing $b''$. 
Since $i_a(a_1,b_1)=1$, we can find $c_1,\ldots,c_h \in \Z$ so that
$$b_1 = [\beta'] + (\sum\nolimits_{i=1}^h c_i a_i) + d [\beta''].$$
See Figure \ref{figure:extendinghomologybasis}.  Let
$\delta$ be the curve indicated there.  Hence
$$[\delta] = [\beta''] - [\alpha_1].$$
Set
$$\beta_1 = T_{\alpha_1}^{c_1+\cdots+c_h+d} T_{\gamma_2}^{-c_2} \cdots T_{\gamma_h}^{-c_h} T_{\delta}^{-d}(\beta').$$
Observe that $i_g(\alpha_i,\beta_1)=\delta_{i1}$.  Also,
\begin{align*}
[\beta_1] = &[\beta'] + (c_1+\cdots+c_h+d)[\alpha_1] + c_2 [\gamma_2] + \cdots + c_h [\gamma_h] + d [\delta] \\
          = &[\beta'] + (c_1+\cdots+c_h+d)[\alpha_1] \\
            &+ c_2 ([\alpha_i] - [\alpha_1]) + \cdots + c_h ([\alpha_h]-[\alpha_1]) + d ([\beta''] - [\alpha_1]) \\
          = &[\beta'] + (c_1 a_1 + \cdots + c_h a_h) + d [\beta''] = b_1,
\end{align*}
as desired.  This reduces us to a previous case, and completes the proof.
\end{proof}

\noindent
Dept. of Mathematics, University of Chicago\\
5734 University Ave.\\
Chicago, Il 60637\\
E-mail: {\tt andyp@math.uchicago.edu}
\medskip


\begin{thebibliography}{}

\bibitem{A}
M. A. Armstrong,
On the fundamental group of an orbit space,
Proc. Camb. Phil. Soc. 61 (1965) 639--646.

\bibitem{Bi1}
J. Birman,
Mapping class groups and their relationship to braid groups,
Comm. Pure Appl. Math. 22 (1969) 213--238.

\bibitem{Bi2}
J. Birman,
On Siegel's modular group,
Math. Ann. 191 (1971) 59--68.

\bibitem{Bi3}
J. Birman,
Braids, Links, and Mapping Class Groups, Annals of Mathematics
Studies, No. 82, Princeton University Press, Princeton, 1974.

\bibitem{Gervais}
S. Gervais,
A finite presentation of the mapping class group of a punctured surface,
Topology 40 (4) (2001) 703--725.

\bibitem{G}
P. Gold,
On the mapping class and symplectic modular group,
thesis, New York University, 1961.

\bibitem{H}
R. Hain,
Torelli groups and Geometry of Moduli Spaces of Curves,
in: C. H. Clemens, J. Kollar (Eds),
Current Topics in Complex Algebraic Geometry,
MSRI publications No. 28,
Cambridge University Press (1995) 97-143.

\bibitem{Ha}
J. L. Harer,
Stability of the homology of the mapping class groups of orientable surfaces,
Ann. of Math. 121 (2) (1985) 215--249.

\bibitem{Hv}
W. J. Harvey,
Geometric structure of surface mapping--class groups,
in: C. T. C. Wall (Ed),
Homological group theory,
London Math. Soc. Lecture Notes Ser., vol. 36,
Cambridge University Press, Cambridge-New York, 1979, pp. 255--269.

\bibitem{HOM}
A. J. Hahn, O. T. O'Meara,
The Classical Groups and K-Theory,
Springer-Verlag, New York, 1989.

\bibitem{I1}
N. V. Ivanov,
Complexes of curves and the Teichm\"{u}ller modular group,
Uspekhi Mat. Nauk 42 (3) (1987) 49--91.

\bibitem{IvanovBook}
N. V. Ivanov,
Subgroups of Teichm\"uller modular groups,
Translations of Mathematical Monographs, No. 115,
AMS, Providence, 1992.

\bibitem{J1}
D. L. Johnson,
A survey of the Torelli Group,
in: Low--dimensional topology (San Francisco, Calif., 1981),
Contemp. Math. 20 (1983) 165--179.

\bibitem{J2}
D. L. Johnson,
The structure of the Torelli group. I. A finite set of
generators for {\sc I},
Ann. of Math. 118 (3) (1983) 423--442.

\bibitem{J3}
D. L. Johnson,
The structure of the Torelli group. II. A characterization of
the group generated by twists on bounding curves,
Topology 24 (2) (1985) 113--126.

\bibitem{J4}
D. L. Johnson,
The structure of the Torelli group. III. The abelianization of
$\Torelli$,
Topology 24 (2) (1985) 127--144.

\bibitem{K}
H. Klingen, 
Charakterisierung der Siegelschen Modulgruppe durch ein endliches System
difinierender Relatimen,
Math. Ann. 144 (1961) 64--82.

\bibitem{Ma}
W. Magnus,
\"{U}ber $n$-dimensionale Gittertransormationen,
Acta. Math. 64 (1934) 353--367.

\bibitem{P}
J. Powell,
Two theorems on the mapping class group of a surface,
Proc. Amer. Math. Soc. 68 (3) (1978) 347--350.

\bibitem{vdB}
B. van den Berg,
On the Abelianization of the Torelli group,
thesis, University of Utrecht, 2003.

\end{thebibliography}
\end{document}